\documentclass[10pt]{article}
\usepackage[margin=1in]{geometry}

\usepackage{amssymb}
\usepackage{amsmath}
\usepackage{amsfonts}
\usepackage{bm}
\usepackage{booktabs}
\usepackage{float}
\usepackage{algorithm}
\usepackage{algorithmic}
\usepackage{multirow}
\usepackage{textcomp}
\usepackage{graphicx}
\usepackage{caption}
\usepackage{subcaption}
\usepackage{microtype}
\usepackage[nocompress]{cite}
\usepackage{hyperref}

\usepackage{amsopn}

\makeatletter

\newcommand{\Rmnum}[1]{\expandafter\@slowromancap\romannumeral #1@}
\makeatother

\begin{document}

\title{A domain decomposition online-learning-enhanced nonlinear elimination preconditioner}
\author{
Pai Zhang\textsuperscript{a}
\qquad
Linyan Gu\textsuperscript{b}
\qquad
Li Luo\textsuperscript{a,*}
\\[1ex]
\small
\textsuperscript{a}Department of Mathematics, University of Macau,
Macao SAR, China.
\\
\small
\textsuperscript{b}School of Mathematics, Sun Yat-sen University,
Guangzhou, China.
\\
\small
\textsuperscript{*}Corresponding author.
E-mail: \texttt{liluo@um.edu.mo}
}

\date{}

\maketitle

\begin{abstract}
Nonlinearly preconditioned inexact Newton methods form an effective class of solvers for large-scale nonlinear algebraic systems arising from the discretization of partial differential equations. A central challenge in nonlinear elimination (NE) preconditioning is the reliable identification of the slowly converging components to be eliminated. Existing selection strategies often rely on problem-specific physical information or user-tuned thresholds applied directly to the raw nonlinear residual, which may contain irregular oscillatory structures near stagnation regions, making the selected bad subset highly sensitive to threshold parameters.
In this work, we propose an online-learning-enhanced NE preconditioner that identifies the bad subset from the dominant structure of the nonlinear residual rather than from the raw residual itself. Residual snapshots are collected online during the stagnation phase of the current Newton solve, and an unsupervised extraction model is trained to capture the principal nonlinear imbalance. We consider both a linear extractor based on principal component analysis and nonlinear extractors based on autoencoder neural networks. Moreover, we integrate the approach into a parallel domain decomposition framework, which trains a local extraction model independently on each subdomain.
The learned residual reconstruction is then used to define the bad subset and guide the nonlinear elimination process. Numerical experiments on  lid-driven cavity flows at Reynolds numbers up to $10{,}000$ show that the proposed method produces more reliable and coherent bad subsets, is robust with respect to both NE and learning parameters, and outperforms the baseline NE preconditioner in terms of the convergence.

\smallskip
\noindent\textbf{Keywords:} Inexact Newton; learning-enhanced nonlinear preconditioning; principal component analysis; neural networks; incompressible Navier--Stokes equations.

\noindent\textbf{MSC 2020:} 49M15; 68T07; 65N55; 76D05.
\end{abstract}

\section{Introduction}
The algebraic systems obtained from the discretization of nonlinear partial differential equations (PDEs) often exhibit strongly localized sources of nonlinear difficulty. In many such problems, the loss of robustness of Newton-type methods is caused by a relatively small set of variables whose nonlinear behavior is substantially more severe than that of the remaining components. These localized components may dominate the residual reduction process and lead to slow convergence or stagnation of the global Newton iteration.

Nonlinear elimination (NE) preconditioning is an effective strategy for improving the robustness and efficiency of Newton-type methods for this class of problems~\cite{Cai2011,Hwang2015,Klawonn2022,Yang2016b,Yang2017,Yang2018,Klawonn2024,Liang2026}. The central idea is to identify the algebraic components that are mainly responsible for the nonlinear difficulty and to treat them separately from the rest of the system. More specifically, the index set of unknowns is divided into a good subset, whose components are sufficiently well balanced and are therefore left to the global Newton iteration, and a bad subset, whose components exhibit strong or slowly converging nonlinear behavior. The residual components associated with the bad subset are approximately eliminated by solving a restricted nonlinear subproblem, and the resulting correction is then incorporated into the outer Newton process.
From this perspective, the NE preconditioner may be viewed as a nonlinear analogue of Gaussian elimination. 

The success of NE preconditioning depends critically on whether the selected bad subset accurately captures the components that slow down the convergence of the original Newton iteration. Many existing selection strategies are based on thresholding the raw nonlinear residual~\cite{Yang2018,Huang2016,Luo2019,Luo2020,Luo2021}. Consequently, the quality of the resulting bad subset can be highly sensitive to user-prescribed parameters and heuristic criteria. An inaccurate selection may either leave the dominant source of nonlinear difficulty untreated or include an unnecessarily large number of degrees of freedom, thereby degrading both the efficiency and robustness of the overall solver. In particular, for complex nonlinear problems, the raw residual often contains oscillatory or highly irregular structures in addition to the dominant nonlinear features. Direct thresholding of such a residual may therefore select many scattered components that are not actually associated with the main source of nonlinear stiffness. Eliminating these isolated components can significantly increase the cost of the restricted nonlinear solve without providing a corresponding improvement in the global convergence. Even worse, the mismatch between the locally updated bad components and the unchanged good components may introduce additional residual jumps across the interface between the two regions. These newly generated jumps may redistribute the nonlinear imbalance, making the subsequent global Newton problem even less balanced.

To address this issue, we propose to exploit the nonlinear residual data generated online during the current stagnating Newton iterations. Instead of selecting the bad subset directly from the raw residual, we use an unsupervised learning model to extract a representation that reveals the dominant nonlinear structures responsible for slow convergence. The extracted dominant structure is then used to define the bad subset and to guide the nonlinear elimination process, so that the preconditioner targets the essential nonlinear imbalance while avoiding many scattered oscillatory components. 

Data-driven extraction of dominant structures from snapshot data has been widely studied in the context of modal decomposition and reduced-order modeling. Classical principal component analysis (PCA), or equivalently proper orthogonal decomposition (POD), identifies energetically dominant coherent modes from a collection of snapshots and has been extensively used to analyze coherent structures in fluid flows~\cite{Sirovich1987,Holmes1996,Taira2017,luo2023pinl} and hyperelastic deformation~\cite{Gong2024}. Related modal-analysis techniques, such as dynamic mode decomposition and spectral POD, further aim to extract dynamically or space-time coherent structures from complex flow data~\cite{Schmid2010,Towne2018}. Beyond linear subspace methods, nonlinear dimensionality-reduction techniques based on neural networks and autoencoders have been developed to learn nonlinear low-dimensional representations of high-dimensional data~\cite{Kramer1991,Hinton2006}. More recently, convolutional autoencoders have been used to extract nonlinear modes and coherent flow structures from fluid-dynamics data~\cite{Murata2020,Fukami2020}, and deep autoencoder manifolds have also been introduced for nonlinear reduced-order modeling~\cite{LeeCarlberg2020}.

Motivated by these developments, we consider both PCA-based and neural-network-based extraction models to identify the dominant structure of the nonlinear residual generated during a stagnating Newton iteration. However, the purpose of the present work differs from that of the aforementioned studies. We do not use the learned representation to construct a reduced-order model, predict system dynamics, or augment a linear preconditioner. Instead, the learning is employed at the nonlinear preconditioning level: the dominant structure extracted from the online residual snapshots replaces the raw residual in defining the bad subset for nonlinear elimination. Moreover, we partition the residual snapshots across the subdomains in a domain decomposition framework \cite{luo2017parallel,Luo2019,Luo2020}, so that the extraction models can be constructed independently on each subdomain, allowing the extraction stage to be performed in parallel. In this way, the NE preconditioner targets the principal nonlinear imbalance while filtering out oscillatory residual components that are less relevant to the slow convergence. Numerical results for high-Reynolds-number lid-driven cavity flows demonstrate that this strategy substantially improves the robustness of NE preconditioning, particularly with respect to the threshold parameter and the stopping tolerance of the inner nonlinear iteration.

The remainder of this paper is organized as follows. 
In Section~\ref{sec:NE}, we review the nonlinear elimination preconditioner and introduce the basic notation used throughout the paper. 
In Section~\ref{sec:algorithms}, we present the proposed DD-OLE-NE-preconditioned inexact Newton algorithm in detail, including online
residual-snapshot collection and the PCA- and autoencoder-based extraction strategies.
 In Section~\ref{sec:experiments}, numerical experiments for the two-dimensional lid-driven cavity flows are presented, including a parametric study of the algorithm and a comparison with the baseline NE preconditioner. Concluding remarks are given in Section~\ref{sec:conclusions}.

\section{The nonlinear elimination preconditioner}
\label{sec:NE}
Consider a nonlinear system of algebraic equations
$F:\mathbb{R}^n\rightarrow\mathbb{R}^n$.
We seek $X^\ast\in\mathbb{R}^n$ such that
\begin{align}\label{nonlinear_system}
F(X^\ast)=0,
\end{align}
starting from an initial guess $X^0\in\mathbb{R}^n$, where
\begin{align*}
F(X)
=
\bigl(
F_1(X),\ldots,F_n(X)
\bigr)^\top,
X
=
\bigl(
X_1,\ldots,X_n
\bigr)^\top.
\end{align*}
We first recall the inexact Newton algorithm with backtracking (IN) \cite{Shadid1997}. Given the current approximate solution $X^k$, a new iterate $X^{k+1}$ is computed via
\begin{align}
X^{k+1}=X^k-\lambda^k S^k,
\end{align}
where the inexact Newton direction $S^k$ satisfies
\begin{align}
\|F(X^k)-F'(X^k)S^k\| \leq \eta^k\|F(X^k)\|.
\end{align}
{Here,} $\eta^k \in [0,1)$ is a forcing term that determines how accurately the Jacobian system needs to be solved.
The step length $\lambda^k\in(0,1]$ is obtained from a standard backtracking line search technique \cite{Dennis1996}. It determines a step size along the inexact Newton direction $S^k$ such that
\begin{align}
f(X^k-\lambda^k S^k)\leq f(X^k) - \alpha \lambda^k \nabla f(X^k)^\top S^k,
\label{linesearch}
\end{align}
where $f=\|F\|^2/2$ is the merit function, and the parameter $\alpha$ is used to ensure that $f$ is reduced sufficiently (herein $\alpha=10^{-4}$).
The nonlinear iteration is stopped if
\begin{align}
\|F(X^k)\| \leq \max\left\{\gamma_a, \gamma_r\|F(X^0)\|\right\}, \label{nonlinearstop}
\end{align}
where $\gamma_a$ and $\gamma_r$ are prescribed absolute and relative tolerances, respectively.

 IN may experience a long stagnation phase when applied to nonlinearly difficult problems. In practice, this slow convergence is often not caused uniformly by all components of the nonlinear residual. Instead, it is typically driven by a relatively small subset of components of $F(X)$ that exhibit stronger nonlinear behavior than the others. These components may dominate the residual reduction process by contributing a large portion of the total nonlinear residual norm, and therefore play a decisive role in the slow convergence or
stagnation of the global Newton iteration.

The central idea of the nonlinear elimination (NE) preconditioner is to identify the components of $F(X)$ that are primarily responsible for the slow convergence of the Newton iteration and to treat them through a restricted nonlinear solve. To introduce the corresponding notation, let
$N=\{1,\ldots,n\}$
denote the index set of all components. The set $N$ is partitioned into a bad subset $N_b$, containing $n_b$ components, and its complementary good subset
$N_g=N\setminus N_b$,
containing $n-n_b$ components. The indices in $N_b$ are associated with variables exhibiting strong or poorly balanced nonlinear behavior, whereas those in $N_g$ correspond to variables with comparatively weak or well-balanced nonlinearities that can be handled effectively by the global Newton iteration.

We define two subspaces
\begin{align*}
&V_b=\left\{v\, | \, v=\left(v_1,\ldots,v_{n}\right)^\top \in \mathbb{R}^n, ~v_{i}=0 \text{ if } i \notin N_b \right\}, \\
&V_g=\left\{v\, | \, v=\left(v_1,\ldots,v_{n}\right)^\top \in \mathbb{R}^n, ~v_{i}=0 \text{ if } i \in N_b \right\}, 
\end{align*}
and the corresponding restriction operators $R_b:\mathbb{R}^n\rightarrow V_b$ and $R_g:\mathbb{R}^n\rightarrow V_g$. Then, $F(X)$ can be partitioned into:
\begin{align}
  F(X)=R_g F(X) + R_b F(X).
  \end{align}

  For a given Newton iterate $X^k$, the NE preconditioner defines a nonlinear mapping $G$ by seeking a locally corrected state $Y = G(X^k)$
such that
\begin{align}
    R_b F(Y)=0 .
    \label{original_precond_system}
\end{align}
The restricted problem~\eqref{original_precond_system} is referred to as the subspace nonlinear system, since only the components corresponding to $N_b$ are targeted for elimination. In practice, this system is solved by an inner Newton iteration.

The NE procedure is embedded within the outer Newton iteration. At each outer step, the current approximation $X^k$ is first used as the initial guess for the subspace solve, and the resulting approximate solution $Y^\ast \approx G(X^k)$ provides a nonlinear correction before the global Newton iteration is continued. Therefore, the subspace nonlinear system does not need to be solved exactly. Starting from
$Y^0 = X^k$,
the inner Newton iteration is terminated once an approximate solution $Y^\ast $ satisfies
\begin{align}
\|R_bF(Y^\ast )\| \leq \gamma_r^s \|R_bF(Y^0)\|,
\end{align}
where $\gamma_r^s$ is the prescribed relative tolerance for the inner nonlinear solve. This inexact treatment reduces the cost of the nonlinear elimination step while still sufficiently damping the bad components that slow down the convergence of the outer Newton iteration.

Although the NE preconditioner has demonstrated considerable success in a variety of applications, its overall effectiveness is still largely influenced by two important factors:
\begin{itemize}
    \item First, the method relies on the reliable identification of the bad subset $N_b$. This selection is highly problem-dependent and remains challenging in practical computations. The reliable construction of $N_b$ is therefore the main challenge addressed in this work. In Section~\ref{sec:algorithms}, we introduce an online-learning-enhanced technique that extracts the dominant nonlinear structures from residual snapshots generated during the current Newton solve. The reconstructed residual is subsequently used in place of the raw residual to define the bad subset, allowing the nonlinear elimination procedure to focus more effectively on the components associated with the principal nonlinear imbalance.
    \item Second, after the nonlinear elimination step is performed, new discontinuities or jumps may be introduced in the residual across the interface between the good and bad regions, or across subdomain boundaries. These newly generated jumps can alter the distribution of nonlinear residuals and may cause previously balanced nonlinear effects to become unbalanced again. To reduce the risk of such interfacial jumps, we employ the modified Newton direction proposed in~\cite{Liu2022}, together with a global line-search procedure. The modified direction accounts for the effect of the local nonlinear correction on the full system, while the line search selects an appropriate step length to ensure sufficient reduction of the global merit function
$f$.
This globalization strategy helps maintain consistency between the locally eliminated components and the remainder of the nonlinear system.
\end{itemize}

 The NE preconditioner using the modified Newton direction is presented in Algorithm~\ref{alg:NE}, stated as a function $X^{new}=\mathrm{NE}(F, X, R_b, R_g, \gamma)$ so that it can later be invoked with a different nonlinear function and different restriction operators (see Algorithm~\ref{alg:nepinl}).
 \begin{algorithm}
  \caption{$X^{new}=\mathrm{NE}(F, X, R_b, R_g, \gamma)$: the NE preconditioning algorithm using a modified Newton direction. The function takes a nonlinear function $F$, an approximate solution $X$, a pair of restriction operators $R_b$ and $R_g$, and a relative tolerance $\gamma$, and returns an updated approximate solution $X^{new}$.}
  \begin{itemize}
    \item[1.] Set $j=0$, start from the initial guess $Y^0=X$, compute $R_bF(Y^0)$ and $\|R_bF(Y^0)\|$.
    \item[2.] While $\|R_bF(Y^j)\| > \gamma\|R_bF(Y^0)\|$, do:
        \begin{itemize}
          \item[(a)] Compute the subspace Jacobian $J^j_b=R_bF'(Y^j)R_b^\top +R_gI_{n\times n}R_g^\top$.
          \item[(b)] Solve $J^j_bS^j_b=R_bF(Y^j)$ for $S^j_b$.
          \item[(c)] Compute the modified global residual $H^j=R_bF'(Y^j-S^j_b)R_b^\top S^j_b+R_gF(Y^j)$.
          \item[(d)] Find the modified Newton direction $\tilde{S}^j$ by approximately solving $F'(Y^j-S^j_b)\tilde{S}^j=H^j$.
          \item[(e)] Update $Y^{j+1}=Y^j-\tilde{\lambda}^j\tilde{S}^j$, where the step length $\tilde{\lambda}^j$ is determined by a backtracking line search technique using the merit function $f=\|F\|^2/2$.
          \item[(f)] Set $j\leftarrow j+1$, compute $R_bF(Y^j)$ and $\|R_bF(Y^j)\|$.
        \end{itemize}
    \item[3.] Return the new approximate solution $X^{new}=Y^j$.
  \end{itemize}
  \label{alg:NE}
\end{algorithm}

\section{A domain decomposition online-learning-enhanced NE preconditioner}\label{sec:algorithms}

\subsection{Motivation}

\begin{figure}[htbp]
    \centering
    \includegraphics[width=0.7\linewidth]{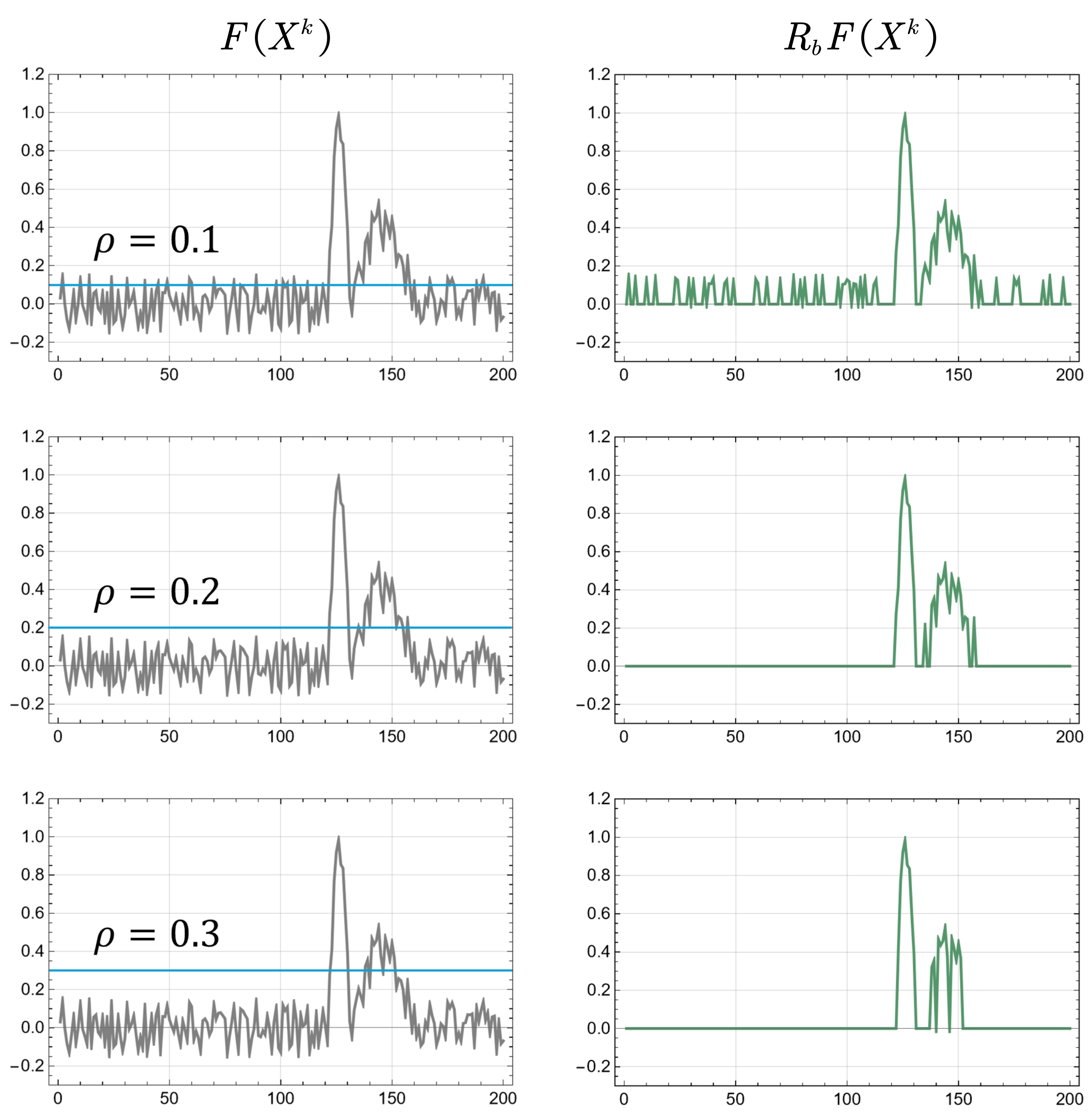}
    \caption{A visual illustration of the residual $F(X^k)$ (left column, with the threshold level $\rho\|F(X^k)\|_\infty$ shown in blue) and the selected bad components $R_b F(X^k)$ (right column) for thresholds $\rho=0.1$, $0.2$ and $0.3$ (top to bottom). A small threshold leaves highly oscillatory protrusions in the bad subset, while a large threshold discards part of the dominant structure.}
    \label{fig:ne_thresholding}
\end{figure}

We illustrate the difficulty of using a classical selection criterion through a simple one-dimensional example. The residual profile $F(X^k)$, shown by the gray curve in the left column of Figure~\ref{fig:ne_thresholding}, contains two pronounced peaks. These peaks represent residual components that remain difficult to reduce even after additional Newton iterations. In many nonlinear problems, such components are associated with localized singular behavior, sharp transitions, or other strongly nonlinear features of the solution. They are therefore regarded as the ``nonlinearly stiff'' components of the residual, which the NE preconditioner aims to identify and eliminate.

An intuitive selection of the bad subset is based on the magnitude of the residual values:
\begin{align}
\label{bad_subset}
N_b :=
\left\{
j\in N
\;\middle|\;
\left|F_j(X^k)\right|
>
\rho\|F(X^k)\|_\infty
\right\},
\end{align}
 where $\rho$ is a user-defined parameter. However, choosing the threshold $\rho$ becomes challenging when the residual contains oscillatory components of relatively small magnitude. Figure~\ref{fig:ne_thresholding} illustrates the bad components $R_bF(X^k)$ selected by different choices of $\rho$. When a small threshold is used, for example $\rho=0.1$, the selected residual contains many disconnected and highly oscillatory protrusions, as shown in the top row. These scattered components enlarge and complicate the resulting subspace nonlinear system, and may make the inner nonlinear solve more difficult rather than more effective.
Increasing the threshold can suppress such protrusions, but this improvement comes at the cost of losing part of the dominant nonlinear structures that should be eliminated. As shown in the case $\rho=0.3$ in the bottom row of Figure~\ref{fig:ne_thresholding}, a larger threshold produces a cleaner bad subset, but may also truncate the main structures associated with the slow convergence. This loss of completeness can weaken the nonlinear elimination effect and eventually affect the robustness of the preconditioner.
These observations suggest that a more effective selection strategy should not be based directly on the raw residual. Instead, it should distinguish the dominant nonlinear structures that need to be eliminated from the highly oscillatory perturbations that contaminate the thresholding process. Such a separation is far from straightforward, however, because the residual $F(X^k)$ inherits the complex spatial and nonlinear features of the underlying nonlinear partial differential equations.

\begin{figure*}[htbp]
    \centering
    \includegraphics[width=\textwidth]{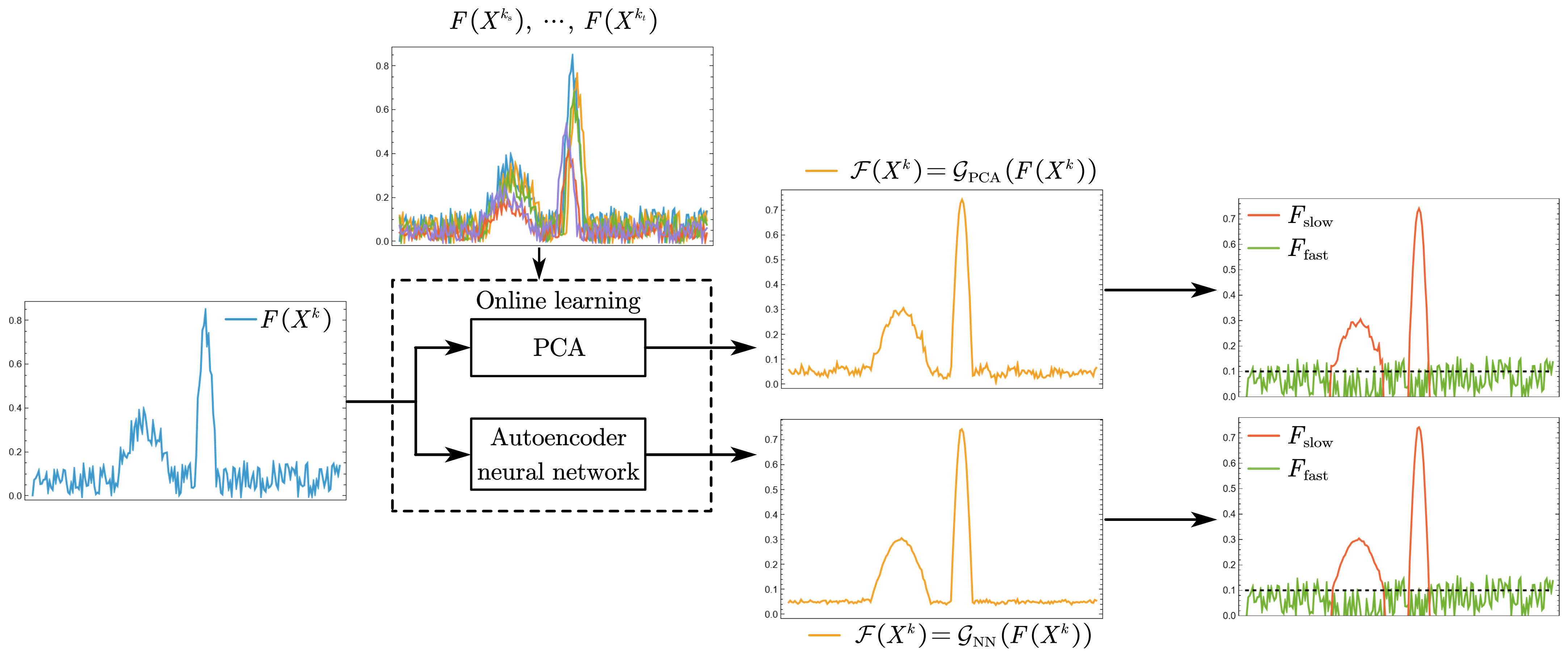}
    \caption{A schematic illustration of the proposed online-learning-based extraction of the dominant residual structure. Nonlinear residual snapshots $F(X^{k_s}),\ldots,F(X^{k_t})$, collected during the stagnation phase, are used to train either a linear PCA-based model or a nonlinear autoencoder.
    The reconstructed residual $\mathcal{F}(X^k)$ preserves the dominant nonlinear structure while filtering out undesirable oscillatory components. Thresholding the reconstruction then identifies the bad components $\mathcal{R}_b\mathcal{F}(X^k)$, shown in red, which are subsequently treated by the NE preconditioner.
}
    \label{fig:arch}
\end{figure*}
 
In this work, we propose a novel strategy to identify the dominant structures associated with localized nonlinearities while filtering out highly oscillatory perturbations in the nonlinear residual. To this end, we introduce a data-driven extraction operator $\mathcal G$ and define the reconstructed residual mapping by
\begin{align}\label{reconstructionF}
\mathcal{F}(X)
=
\mathcal{G}(F(X))
=
\bigl(
\mathcal{F}_1(X),\ldots,\mathcal{F}_n(X)
\bigr)^\top
\in\mathbb{R}^n,
\end{align}
where $\mathcal{G}(\cdot)$ is a data-driven extraction operator learned online from the collected residual snapshots. We consider two realizations: a linear model based on principal component analysis and nonlinear models based on neural networks.

Applying $\mathcal{G}$ to the residual produces the reconstructed residual $\mathcal{F}(X^k)$, whose dominant nonlinear structures are preserved while oscillatory perturbations are suppressed. Figure~\ref{fig:arch} illustrates the proposed online-learning-based extraction process. Compared with the raw residual $F(X^k)$, the reconstructed residual $\mathcal{F}(X^k)$ more clearly exposes the dominant structures associated with localized nonlinearities while removing high-frequency oscillations, thereby allowing a more reliable identification of the bad subset. Replacing the raw residual with its reconstruction in the selection criterion yields
\begin{align}
\label{new_bad_subset}
\mathcal{N}_b :=
\left\{
j\in N
\;\middle|\;
\left|\mathcal{F}_j(X^k)\right|
>
\rho\|\mathcal{F}(X^k)\|_\infty
\right\}.
\end{align}

As illustrated in Figure~\ref{fig:arch}, thresholding the reconstructed residual produces the bad components
$\mathcal{R}_b\mathcal{F}(X^k)$, which capture the dominant nonlinear structures while avoiding the scattered discontinuous protrusions generated by thresholding the raw residual (top row of Figure~\ref{fig:ne_thresholding}). Consequently, the bad-subset selection is expected to become less sensitive to the threshold parameter, thereby improving the robustness of the NE preconditioner.

With the identified bad subset $\mathcal{N}_b$, we define the corresponding restriction operators $\mathcal{R}_b$ and $\mathcal{R}_g$, analogously to $R_b$ and $R_g$, and decompose the nonlinear residual as
\begin{align}\label{decomposition}
F = F_{\mathrm{slow}} + F_{\mathrm{fast}}.
\end{align}
Here,
\begin{align*}
F_{\mathrm{slow}}
:=
\mathcal{R}_b\mathcal{F}
\end{align*}
represents the dominant residual structures associated with localized nonlinearities that are primarily responsible for the slow convergence of the global Newton iteration. Throughout this paper, these components are referred to as the \emph{slow components}. The remaining part,
\begin{align*}
F_{\mathrm{fast}}
:=
F-\mathcal{R}_b\mathcal{F},
\end{align*}
contains residual components whose nonlinearities are comparatively well balanced and can therefore be handled efficiently by the global Newton iteration.

To reduce the slow components, we introduce the following modified subspace nonlinear system:
\begin{align}\label{NN_precond_system}
\mathcal{R}_b\mathcal{F}(Y)=0.
\end{align}
Compared with the original subspace system~(\ref{original_precond_system}), the modified system is generally easier to solve because the reconstructed residual removes many oscillatory components that unnecessarily complicate the inner nonlinear iteration. As a result, the nonlinear elimination preconditioner becomes more robust while preserving its ability to target the dominant sources of nonlinear difficulty.

In practice, incorporating a learned extraction operator into nonlinear preconditioning requires resolving several key issues:
\begin{enumerate}
    \item How to collect the training data?
    \item How to efficiently train the extraction models on large-scale residual data in parallel?
    \item How to define the extraction models to optimally extract the dominant structure of the residual?
    \item How to integrate the extraction models into the NE preconditioner in Algorithm~\ref{alg:NE}?
\end{enumerate}
We address these questions in the subsequent sections.

\subsection{Collection of training data}
\begin{figure}[htbp]
  \centering
    \includegraphics[width=0.65\linewidth]{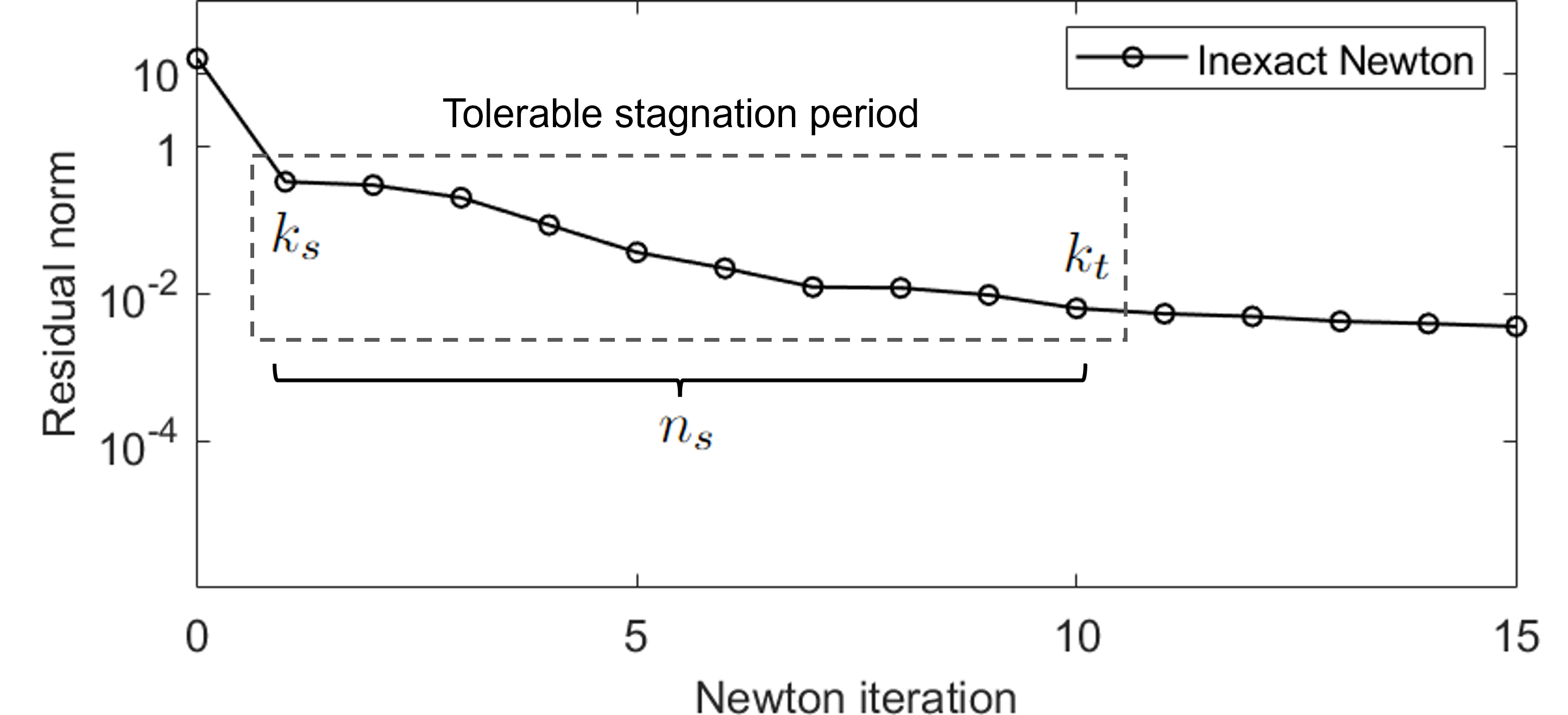}
  \caption{A typical residual curve obtained using IN when the method fails to converge. The part in the dashed box is considered as the tolerable stagnation period $[k_{s}, k_{t}]$ with length $n_{s}=k_{t}-k_{s}+1$.}
  \label{fig:data_collection}
\end{figure}

Nonlinear preconditioning is only necessary when the convergence of IN is problematic. Figure \ref{fig:data_collection} shows a typical residual curve when the method becomes trapped in a state of stagnation. We consider the stagnation to be tolerable until reaching a step $k_{t}$ where the nonlinear preconditioner is needed to accelerate the convergence. Letting $k_{s}$ denotes the starting step of the stagnation, we define a tolerable stagnation period $[k_{s}, k_{t}]$ with length $n_{s}=k_{t}-k_{s}+1$. During this period, there often exist dominant coherent structures within the set of nonlinear residual vectors generated by the Newton iterations:
  \begin{align}
   S_F = \left\{ F(X^{k_{s}}), F(X^{k_{s} + 1}), \ldots, F(X^{k_{t}}) \right\}.
  \end{align}
 These residual vectors offer useful hidden information that is associated with the slow convergence of IN. We collect the set $S_F$ as a training dataset for the proposed unsupervised learning approach to identify the hidden information. 

The size of the dataset is determined by the length of the tolerable stagnation period $n_s$. In practice, $k_s$ is chosen as the first iteration at which the decrease of the residual norm becomes noticeably slow, while $k_t$ is selected after the stagnation has persisted long enough to justify the use of nonlinear preconditioning.

\subsection{Parallel training based on domain decomposition}
\label{subsec:dd}

\begin{figure}[htbp]
  \centering
    \includegraphics[width=0.75\linewidth]{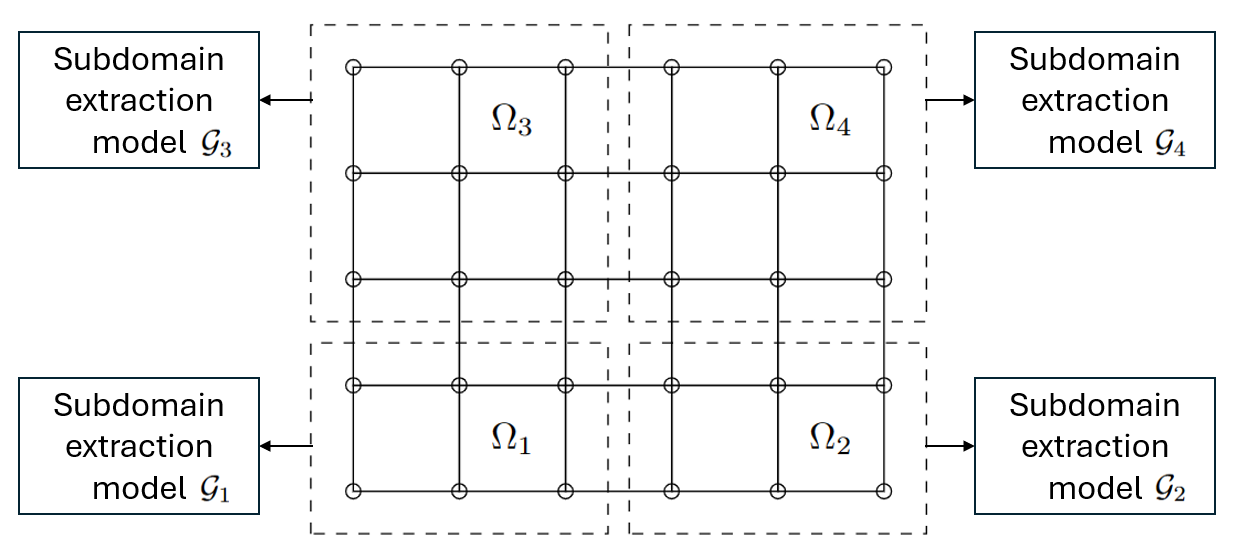}
  \caption{Subdomain-wise extraction models used in
the proposed DD-OLE-NE preconditioner. The computational domain is
decomposed into four nonoverlapping subdomains
$\Omega_1,\dots,\Omega_4$. On each subdomain
$\Omega_i$, local residual snapshots are used to construct an
independent extraction model $\mathcal{G}_i$.}
  \label{fig:subdomains}
\end{figure}

To enable efficient online training for large-scale nonlinear systems, we adopt a domain-decomposition-based learning strategy. As shown in Figure~\ref{fig:subdomains}, the computational domain $\Omega$ is partitioned into $n_p$ nonoverlapping subdomains,
\begin{align}
    \Omega=\bigcup_{i=1}^{n_p}\Omega_i,
    \qquad
    \Omega_i\cap\Omega_j=\emptyset,
    \quad i\neq j.
\end{align}
Accordingly, the global solution vector $X\in\mathbb{R}^n$ is partitioned into subdomain vectors. Let $R_i^0$ denote the restriction operator from the global vector to the degrees of freedom associated with the nonoverlapping subdomain $\Omega_i$. The local solution vector on $\Omega_i$ is therefore defined as
$X_i = R_i^0 X$.
After a suitable ordering of the global degrees of freedom, the restriction can be written formally as
\begin{align}
    X_i
    =
    R_i^0 X
    =
    (I \quad 0)
    \begin{pmatrix}
        X_i \\
        X\backslash X_i
    \end{pmatrix},
\end{align}
where $X\backslash X_i$ denotes the components of $X$ outside the subdomain $\Omega_i$. Let $n_i$ denote the number of degrees of freedom associated with
$\Omega_i$.

The same decomposition is applied to the nonlinear residual snapshots collected during the stagnation period. Instead of training a single global extraction model using the full residual dataset, the snapshots are restricted to each subdomain to form the local training dataset
\begin{align}
    S_{F,i}
    =
    \left\{
    R_i^0F(X^{k_s}),
    R_i^0F(X^{k_s+1}),
    \ldots,
    R_i^0F(X^{k_t})
    \right\},
\end{align}
where $i=1,\ldots,n_p$.
An independent extraction model $\mathcal{G}_i$ is then constructed and trained using $S_{F,i}$ on the corresponding processor
\begin{align}
    \mathcal{G}_i:\mathbb{R}^{n_i}\rightarrow\mathbb{R}^{n_i},
\end{align}
where $\mathcal{G}_i$ may be realized either by a local PCA model or by a
local autoencoder neural network.

In the proposed DD-OLE-NE preconditioner, for each residual snapshot
$F(X^k)$, $k=k_s,\ldots,k_t$, its restriction to the $i$-th
nonoverlapping subdomain is $R_i^0F(X^k)$. The corresponding local
reconstruction is
$\mathcal{G}_i\bigl(R_i^0F(X^k)\bigr)$.
The global reconstructed residual is assembled as
\begin{align}
\mathcal{F}(X^k)
=
\sum_{i=1}^{n_p}
(R_i^0)^\top
\mathcal{G}_i\bigl(R_i^0F(X^k)\bigr),
\end{align}
where $(R_i^0)^\top$ injects the local reconstruction into the global
vector.

The use of domain decomposition in the learning-based extraction step is
motivated by several considerations.

\begin{itemize}
    \item[1)] The nonlinear residual is often location-sensitive: the
    dominant structures responsible for Newton stagnation may appear only
    in a few localized regions. Applying a global PCA or a global
    autoencoder to the entire residual may mix information from unrelated
    regions and obscure such localized structures. In the spirit of the
    domain decomposition PCA strategy for location-sensitive data
    discussed in~\cite{LiCai2021}, local extraction models can better
    preserve spatially localized features.

    \item[2)] The dimension of each local residual snapshot
    $R_i^0F(X^k)$ is much smaller than that of the global residual
    $F(X^k)$. Therefore, the local PCA or neural network training
    problems are cheaper in both memory and computational cost. This is
    particularly important when the nonlinear system is large and the
    residual snapshots are collected online within a single Newton solve.

    \item[3)] The local extraction models are independent across
    subdomains. Each processor can train its own model
    $\mathcal{G}_i$ using only the local dataset $S_{F,i}$, and no
    global training problem is required \cite{Jagtap2020cPINN,Jagtap2020XPINN,Shukla2021}. Thus, the learning stage fits
    naturally into the domain decomposition parallel framework already
    used by the nonlinear solver.
\end{itemize}

After the global reconstructed residual $\mathcal{F}(X)$ is assembled,
the bad subset is identified by using~(\ref{new_bad_subset}).
In this way, the proposed method combines the scalability of domain
decomposition with the structure-identification capability of PCA or
neural network based models.

\subsection{Extraction models}\label{subsec:recon_models}

For the extraction of the dominant structure of the residual vectors, we consider two classes of models: a linear model based on PCA, and nonlinear models based on neural networks, for which we use a multilayer perceptron and a convolutional neural network. 

\subsubsection{Principal component analysis based model}\label{subsubsec:pca}

In this method, PCA is applied independently on each subdomain to
construct a local linear extraction operator from the residual snapshots
collected during the stagnation phase. For each $i=1,\ldots,n_p$, the residual snapshots are
restricted to $\Omega_i$ and assembled into the local snapshot matrix $\mathbf{F}_i \in \mathbb{R}^{n_i\times n_s}$, where
\begin{align*}
  \mathbf{F}_i
  =
  \left[
  R_i^0F(X^{k_s}),\,
  R_i^0F(X^{k_s+1}),\,
  \ldots,\,
  R_i^0F(X^{k_t})
  \right],
\end{align*}
and $n_s=k_t-k_s+1$ is the number of residual snapshots collected
during the tolerable stagnation period.

For each subdomain, PCA seeks an matrix with orthonormal columns
$P_i\in\mathbb{R}^{n_i\times d}$, satisfying
$P_i^\top P_i=I_d$, where $d\ll n_i$ is the prescribed number of
retained modes. The columns of $P_i$ span a low-dimensional local
subspace that retains the dominant variability of the residual
snapshots $R_i^0F(X^k)$, $k=k_s,\ldots,k_t$.

The mean of the local residual snapshots is defined by
\begin{align*}
  \bar{F}_i
  =
  \frac{1}{n_s}
  \sum_{\ell=k_s}^{k_t} R_i^0F(X^\ell)
  \in\mathbb{R}^{n_i}.
\end{align*}
For $k=k_s,\ldots,k_t$, the corresponding centered local residual
vector is
\begin{align*}
  \hat{F}_i^k
  =
  R_i^0F(X^k)-\bar{F}_i,
\end{align*}
and the centered local snapshot matrix is
\begin{align*}
  \widehat{\mathbf{F}}_i
  =
  \left[
  \hat{F}_i^{k_s},\,
  \hat{F}_i^{k_s+1},\,
  \ldots,\,
  \hat{F}_i^{k_t}
  \right]
  \in\mathbb{R}^{n_i\times n_s}.
\end{align*}

The local PCA basis is obtained from the singular value decomposition
(SVD)
\begin{align*}
  \widehat{\mathbf{F}}_i
  =
  \hat{U}_{i}\,
  \hat{\Sigma}_{i}\,
  \hat{V}_{i}^\top,
  \qquad i=1,\ldots,n_p,
\end{align*}
where
$\hat{U}_{i}\in\mathbb{R}^{n_i\times n_i}$ and
$\hat{V}_{i}\in\mathbb{R}^{n_s\times n_s}$ are orthogonal matrices,
and
$\hat{\Sigma}_{i}\in\mathbb{R}^{n_i\times n_s}$ contains the
singular values
\begin{align*}
  \sigma_{i}^{1}
  \ge
  \sigma_{i}^{2}
  \ge
  \cdots
  \ge
  \sigma_{i}^{\min(n_i,n_s)}
  \ge 0
\end{align*}
in nonincreasing order. The local projection matrix $P_i$ is formed
from the first $d$ columns of $\hat{U}_{i}$, with
$d\leq\min(n_i,n_s-1)$.

The leading left singular vectors represent the most energetic
structures shared by the local residual snapshots. Retaining only the
first $d$ modes therefore preserves the dominant coherent structure
of the local residual while suppressing low-energy oscillatory
perturbations. Since the number of snapshots $n_s$ is typically small
and the SVD is performed independently on each subdomain, the local PCA
construction has a modest computational cost and can be carried out in
parallel.

Using the local PCA bases, the global extracted residual
structure is obtained by
\begin{align}\label{pca_operator}
\mathcal{F}(X^k)
&=
\mathcal{G}_{\mathrm{PCA}}\bigl(F(X^k)\bigr)\nonumber\\
&=
\sum_{i=1}^{n_p}
(R_i^0)^\top
\mathcal{G}_{i,\mathrm{PCA}}
\bigl(R_i^0F(X^k)\bigr)
\end{align}
where the local PCA reconstruction on subdomain $i$ is given by
\begin{align}\label{local_pca_operator}
  \mathcal{G}_{i,\mathrm{PCA}}\bigl(R_i^0F(X)\bigr)
  :=
  P_iP_i^\top
  \left(
  R_i^0F(X)-\bar{F}_i
  \right)
  +
  \bar{F}_i .
\end{align}
The
number of retained modes $d$ is the only model parameter, and the same
PCA-based reconstruction is carried out independently on all subdomains.

\subsubsection{Multilayer perceptron based model}
The multilayer perceptron (MLP) is a feed-forward fully connected network that applies linear and nonlinear transformations to the input recursively. 

For the $i$-th nonoverlapping subdomain, this model constructs an MLP-based
autoencoder to learn a local nonlinear extraction operator that maps the
local residual snapshot to its reconstructed dominant residual structure.
More precisely, the local MLP-based extraction operator is defined by
\begin{align}
    \mathcal{G}_{i,\mathrm{MLP}}(\,\cdot\,;\theta_i)
    =
    \mathcal{D}_{i}(
        \mathcal{E}_{i}(\,\cdot\,;\theta_i^E);
        \theta_i^D
    ),
\end{align}
where $\mathcal{E}_{i}$ and $\mathcal{D}_{i}$ denote the encoder and
decoder on the $i$-th subdomain, respectively. The corresponding global
MLP-based extraction operator is defined by assembling all local
reconstructions:
$$
    \mathcal{G}_{\mathrm{MLP}}(F(X^k);\Theta)
    =
    \sum_{i=1}^{n_p}
    (R_i^0)^\top
    \mathcal{G}_{i,\mathrm{MLP}}
    \bigl(R_i^0F(X^k);\theta_i\bigr),
$$
where $\Theta=\{\theta_i\}_{i=1}^{n_p}$
collects the trainable parameters over all subdomains.

Specifically, the input of the local autoencoder on the $i$-th
subdomain is the restricted residual snapshot
$$
    z_i^{0,k}
    =
    R_i^0F(X^k),
    \qquad k=k_s,\ldots,k_t .
$$
The encoder maps this local residual vector to a low-dimensional
bottleneck representation. For $\ell=1,\ldots,L_E$, we define
$$
    z_i^{\ell,k}
    =
    \sigma\left(
        W_i^\ell z_i^{\ell-1,k}+b_i^\ell
    \right),
$$
where $L_E$ is the number of encoder layers. $W_i^\ell$ denotes the weight,
$b_i^\ell$ is the bias, and $\sigma(\cdot)$ is a nonlinear activation
function. In our implementation, $\sigma$ in MLP is set to be ReLU.
The bottleneck variable is
then given by
$$
    c_i^k
    =
    z_i^{L_E,k}
    \in \mathbb{R}^{r_i},
    \qquad
    r_i \ll n_i.
$$
This bottleneck representation is intended to retain the
dominant nonlinear structure of the local residual while suppressing
oscillatory perturbations.

The decoder maps the bottleneck variable back to the original local
residual space. For $\ell=L_E+1,\ldots,L-1$, we set
$$
    z_i^{\ell,k}
    =
    \sigma\left(
        W_i^\ell z_i^{\ell-1,k}+b_i^\ell
    \right),
$$
and the output layer gives the reconstructed local residual structure
$$
\mathcal{G}_{i,\mathrm{MLP}}
\bigl(R_i^0F(X^k);\theta_i\bigr)
=
W_i^L z_i^{L-1,k}+b_i^L .
$$
Equivalently, the local reconstruction can be written in the
autoencoder form
$$
\mathcal{G}_{i,\mathrm{MLP}}
\bigl(R_i^0F(X^k);\theta_i\bigr)
=
\mathcal{D}_{i}
\left(
    \mathcal{E}_{i}
    \bigl(R_i^0F(X^k);\theta_i^E\bigr);
    \theta_i^D
\right).
$$
Here,
$$
    \theta_i^E
    =
    \{W_i^\ell,b_i^\ell\}_{\ell=1}^{L_E},
    \qquad
    \theta_i^D
    =
    \{W_i^\ell,b_i^\ell\}_{\ell=L_E+1}^{L},
$$
and
$\theta_i
    =
    \{\theta_i^E,\theta_i^D\}
    =
    \{W_i^\ell,b_i^\ell\}_{\ell=1}^{L}
$
denotes all trainable parameters of the local MLP autoencoder.

The local MLP autoencoder is trained in an unsupervised manner by
minimizing the reconstruction error over the residual snapshots on the
$i$-th subdomain:
\begin{align}\label{optimize_mlp}
    \min_{\theta_i}
    \frac{1}{n_s}
    \sum_{k=k_s}^{k_t}
    \left\|
        \mathcal{G}_{i,\mathrm{MLP}}
        \bigl(R_i^0F(X^k);\theta_i\bigr)
        -
        R_i^0F(X^k)
    \right\|_2^2.
\end{align}
After the
local autoencoders have been trained independently on all subdomains, the
global extracted residual structure is obtained by
\begin{align}
    \mathcal{F}(X^k)
    =
    \mathcal{G}_{\mathrm{MLP}}
    \bigl(F(X^k);\Theta^\ast\bigr),
\end{align}
where
$\Theta^\ast=\{\theta_i^\ast\}_{i=1}^{n_p}$ collects the minimizers of (\ref{optimize_mlp}).
The reconstructed residual $\mathcal{F}(X^k)$ is then used in place of
the raw residual $F(X^k)$ to identify the bad subset for nonlinear
elimination.

\subsubsection{Convolutional neural network based model}
We also consider a CNN-based autoencoder to define a nonlinear extraction
operator on each subdomain. Unlike the MLP autoencoder,
which treats the local residual snapshot as an unstructured vector, the
CNN autoencoder exploits the spatial arrangement of the degrees of freedom
on the underlying mesh. This is particularly useful when the dominant
nonlinear structures of the residual are localized and spatially coherent.

Before applying the CNN autoencoder, the
local residual vector is reshaped into a tensor form
$$
Z_i^{0,k}
=
\mathcal{T}_i\bigl(R_i^0F(X^k)\bigr)
\in\mathbb{R}^{h_i\times w_i\times c}.
$$
where $\mathcal{T}_i$ denotes the reshaping operator from the local
residual vector to the spatial tensor on $\Omega_i$, $h_i$ and $w_i$
are the local mesh dimensions, and $c$ is the number of physical
components or fields.

The local CNN-based extraction operator is defined by an encoder-decoder
architecture,
$$
    \mathcal{G}_{i,\mathrm{CNN}}(\cdot;\theta_i)
    =
    \mathcal{T}_i^{-1}
    \circ
    \mathcal{D}_i
    \left(
        \mathcal{E}_i
        \left(
            \mathcal{T}_i(\cdot);
            \theta_i^E
        \right);
        \theta_i^D
    \right),
$$
where $\mathcal{E}_i$ and $\mathcal{D}_i$ denote the convolutional
encoder and decoder on the $i$-th subdomain, respectively, and
$\theta_i=\{\theta_i^E,\theta_i^D\}$
collects the trainable parameters.

The encoder maps the input tensor $Z_i^{0,k}$ to a spatially coarsened 
bottleneck representation through a sequence of convolutional layers. For
$\ell=1,\ldots,L_E$, 
$$
    Z_i^{\ell,k}
    =
    \sigma
    \left(
        \mathrm{Conv}_i^\ell
        \left(
            Z_i^{\ell-1,k};
            K_i^\ell
        \right)
        +
        b_i^\ell
    \right),
$$
where $\mathrm{Conv}_i^\ell$ denotes the $\ell$-th convolutional
operator on the $i$-th subdomain, $K_i^\ell$ is the convolution kernel,
$b_i^\ell$ is the bias, and $\sigma(\cdot)$ is a nonlinear activation
function. In our implementation, each convolutional layer uses a stride of two, so
that the spatial resolution is approximately halved in each direction per
encoder layer, and $\sigma$ is taken as the Mish activation.
The bottleneck variable is given by
$$
    C_i^k
    =
    Z_i^{L_E,k}
    \in
    \mathbb{R}^{h_i^c\times w_i^c\times c_i^c},
    \qquad
    h_i^c < h_i, \quad w_i^c < w_i .
$$
This bottleneck representation is designed to retain the dominant
spatially coherent nonlinear structure of the local residual while
suppressing oscillatory perturbations.

The decoder then maps the bottleneck representation back to the original
local residual space. For $\ell=L_E+1,\ldots,L-1$, we define
$$
    Z_i^{\ell,k}
    =
    \sigma
    \left(
        \mathrm{ConvT}_i^\ell
        \left(
            Z_i^{\ell-1,k};
            K_i^\ell
        \right)
        +
        b_i^\ell
    \right),
$$
where $\mathrm{ConvT}_i^\ell$ denotes a transposed convolutional
operator. The transposed convolutions also use a stride of two, with the output
padding chosen according to the parity of $h_i$ and $w_i$ so that the
reconstructed tensor recovers the original spatial dimensions.
The output layer produces the reconstructed residual tensor
$$
    \widehat{Z}_i^k
    =
    \mathrm{ConvT}_i^L
    \left(
        Z_i^{L-1,k};
        K_i^L
    \right)
    +
    b_i^L .
$$
The reconstructed local residual structure is then obtained by reshaping
the tensor back to vector form:
$$
\mathcal{G}_{i,\mathrm{CNN}}
\bigl(R_i^0F(X^k);\theta_i\bigr)
=
\mathcal{T}_i^{-1}
\left(
    \widehat{Z}_i^k
\right).
$$

The corresponding global CNN-based extraction operator is defined by
assembling the local reconstructed residual structures:
$$
    \mathcal{G}_{\mathrm{CNN}}
    \bigl(F(X^k);\Theta\bigr)
    =
    \sum_{i=1}^{n_p}
    (R_i^0)^\top
    \mathcal{G}_{i,\mathrm{CNN}}
    \bigl(R_i^0F(X^k);\theta_i\bigr),
$$
where $\Theta=\{\theta_i\}_{i=1}^{n_p}$.

The local CNN autoencoder is trained in an unsupervised manner similar to (\ref{optimize_mlp}).
The global extracted residual structure is obtained by
$$
    \mathcal{F}(X^k)
    =
    \mathcal{G}_{\mathrm{CNN}}
    \bigl(F(X^k);\Theta^\ast\bigr),
$$
with
$\Theta^\ast=\{\theta_i^\ast\}_{i=1}^{n_p}$ denoting the collection of trained parameters.
The reconstructed residual $\mathcal{F}(X^k)$ is then used in place of
the raw residual $F(X^k)$ to identify the bad subset for nonlinear
elimination.

In the implementation, each residual component is standardized before being passed to the MLP or CNN, and the reconstructed output is inverse-transformed before the global residual reconstruction is assembled. Moreover, the extraction models are used only to identify the bad subset
and to evaluate the residual vectors in the nonlinear elimination phase.
That is, when Algorithm~\ref{alg:NE} is invoked with $\mathcal{F}$ in
Phase~3, the residual evaluations are based on the reconstructed residual. In contrast, the
Jacobian matrices in steps (a), (c), and (d) are assembled from the
original nonlinear residual $F$, restricted to the selected bad subset.
Equivalently, the Jacobian $F'$ of the original residual is used as an approximation to $\mathcal{F}'$. 
Hence, no differentiation of the extraction models is required, and $F'$ serves as an approximation to $\mathcal{F}'$, which is justified since the reconstruction preserves the dominant structure of the residual.

The complete domain decomposition online-learning-enhanced nonlinear elimination (DD-OLE-NE) preconditioned inexact Newton method is summarized in Algorithm~\ref{alg:nepinl}.

 \begin{algorithm*}
\caption{The domain decomposition online-learning-enhanced nonlinear
elimination (DD-OLE-NE) preconditioned inexact Newton method. Given an
initial guess $X^0$, the nonlinear operator $F$, the nonoverlapping
subdomain restriction operators $R_i^0$, $i=1,\ldots,n_p$, the absolute
and relative tolerances $\gamma_a$ and $\gamma_r$, the relative tolerance
$\gamma_r^s$ for the inner nonlinear iteration, the residual collection
window $[k_s,k_{t}]$, the restart step $k_{re}\in[k_s,k_{t}]$, and the bad subset threshold $\rho$.}
\label{alg:nepinl}
\begin{description}

\item[Initialization.]
Set $k=0$, compute $F(X^0)$ and $\|F(X^0)\|$. Initialize the local
residual snapshot datasets $S_{F,i}=\emptyset,  i=1,\ldots,n_p$.

\item[Phase 1.] The global IN phase:\\
While $\|F(X^k)\| > \max\left\{\gamma_a,\gamma_r\|F(X^0)\|\right\}$, do:
\begin{itemize}
    \item[(1.1)] If $k_s\leq k\leq k_{t}$, collect the local nonlinear
    residual snapshots $R_i^0F(X^k)$ and update
    ${S}_{F,i} =
        {S}_{F,i}
        \cup
        \left\{
        R_i^0F(X^k)
        \right\}$.
    If $k=k_{re}$, store $X^{k_{re}}$.
    If $k=k_{t}$, go to Phase 2.
    \item[(1.2)] Inexactly solve the global Newton system
    $F'(X^k)S^k=F(X^k)$.
    \item[(1.3)] Compute the step length $\lambda^k$ using the cubic
    backtracking line search.
    \item[(1.4)] Update $ X^{k+1}=X^k-\lambda^kS^k$.
    \item[(1.5)] Set $k=k+1$, compute $F(X^k)$ and $\|F(X^k)\|$.
\end{itemize}
\item[Phase 2.] The online training phase:
\begin{itemize}
    \item[(2.1)] For each subdomain $\Omega_i$, construct and train a
    local extraction model $\mathcal{G}_i
        \in
        \left\{
        \mathcal{G}_{i,\mathrm{PCA}},
        \mathcal{G}_{i,\mathrm{MLP}},
        \mathcal{G}_{i,\mathrm{CNN}}
        \right\}$ from the collected residual snapshots.
\item[(2.2)] Define the
    reconstructed global residual by assembling the local reconstructed
    residual structures:
    $\mathcal{F}(X)
        =
        \sum_{i=1}^{n_p}
        (R_i^0)^\top
        \mathcal{G}_i
        \left(
        R_i^0F(X)
        \right)$.

    \item[(2.3)] Store the trained local extraction models
    $\mathcal{G}_i(\,\cdot\,)$, $i=1,\ldots,n_p$, and go to
    Phase 3.

\end{itemize}

\item[Phase 3.] The nonlinear elimination phase:
\begin{itemize}

    \item[(3.1)] 
    Set $Y^0=X^{k_{re}}$,
    and compute the reconstructed residual
    $\mathcal{F}(Y^0)
        =
        \sum_{i=1}^{n_p}
        (R_i^0)^\top
        \mathcal{G}_i
        \left(
        R_i^0F(Y^0)
        \right)$.

    \item[(3.2)] Identify the bad subset from the reconstructed residual by
        $\mathcal{N}_b
        =
        \left\{
        j\in N
        \;\middle|\;
        \left|\mathcal{F}_j(Y^0)\right|
        >
        \rho
        \|\mathcal{F}(Y^0)\|_\infty
        \right\}$,
        and construct the corresponding restriction operators
        $\mathcal{R}_b$ and $\mathcal{R}_g$.
    \item[(3.3)] Approximately solve the reconstructed subspace nonlinear
    system $\mathcal{R}_b\mathcal{F}(Y^\ast)=0$
    by calling Algorithm~\ref{alg:NE},
    $Y^\ast
        =
        \mathrm{NE}
        \left(
        \mathcal{F},
        Y^0,
        \mathcal{R}_b,
        \mathcal{R}_g,
        \gamma_r^s
        \right)$.

    \item[(3.4)] Set $k\leftarrow k_t+1$ and the new approximate solution
    $X^k=Y^\ast$, compute $F(X^k)$ and $\|F(X^k)\|$, and restart from Phase 1.
\end{itemize}
\end{description}
\end{algorithm*}

\section{Numerical experiments}\label{sec:experiments}
We evaluate the proposed method on the two-dimensional lid-driven cavity flow at high Reynolds numbers.
Let $\Omega=(0, 1)\times(0, 1)$. The problem can be modeled by the Navier-Stokes equations in the velocity-vorticity formulation:
\begin{equation}
  \left\{\begin{array}{llll}\label{Navier-Stokes}
            \displaystyle  -\Delta u - \frac{\partial \omega}{\partial y}&=0, \quad &\mbox{in}\quad\Omega,\\[5pt]
            \displaystyle  -\Delta v + \frac{\partial \omega}{\partial x}&=0, \quad &\mbox{in}\quad\Omega,\\[5pt]
            \displaystyle -\frac{1}{Re}\Delta \omega + u\frac{\partial \omega}{\partial x} + v\frac{\partial \omega}{\partial y}&=0,  \quad &\mbox{in}\quad\Omega,
         \end{array}\right.
\end{equation}
where $u$ and $v$ are the velocity fields in the $x$- and $y$-directions, respectively, and
\begin{align}
  \displaystyle \omega = \frac{\partial v}{\partial x} - \frac{\partial u}{\partial y}
\end{align}
is the vorticity normal to the $xy$-plane. 
The top boundary $\Gamma_{lid}$ represents a lid moving with velocity $u=1$ in the positive $x$-direction. We impose a no-slip and no-penetration boundary condition on all walls,
\begin{equation}
  \left\{\begin{array}{llll}\label{BC_driven_cavity}
    u = 1, \quad &\mbox{on}\quad \Gamma_{lid},\\[5pt]
    u = 0, \quad &\mbox{on}\quad \partial \Omega\setminus\Gamma_{lid},\\[5pt]
    v =0,  \quad &\mbox{on}\quad \partial \Omega,\\[5pt]
    \omega = \frac{\partial v}{\partial x} - \frac{\partial u}{\partial y}, \quad &\mbox{on}\quad \partial \Omega.
\end{array}\right.
\end{equation}
The boundary condition for the vorticity is discretized by a second-order approximation that uses the mesh points neighboring the boundary~\cite{Prudencio2006}.

The Reynolds number $Re$  quantifies the relative importance of inertial forces to viscous forces. 
A standard central second-order finite difference scheme is used for the discretization of both the Laplacian operators and the first order partial derivatives in (\ref{Navier-Stokes}). Let $\Omega$ be covered by a $M\times M$ mesh, then each point $p_{ij}=(x_i, y_j)$ is located at the position $x_i=(i-1)h$, $y_j=(j-1)h$ with $h=1/(M-1)$ and  $i,j=1,\ldots,M$. In this work, we consider the point-block ordering to build up the large sparse nonlinear system of algebraic equations (\ref{nonlinear_system}), in which the unknown variables $u_{ij}$, $v_{ij}$, $\omega_{ij}$ associated with a mesh point $p_{ij}$ are always together in a $3\times 3$ block. For this multi-variable system, the bad subset is defined by a pointwise approach \cite{Luo2020} in the NE preconditioner. 

The numerical experiments are carried out on a computer with an Intel Xeon 6248 2.50GHz CPU and eight NVIDIA Tesla V100-SXM2 32GB GPUs.
A zero vector is used as the initial guess, i.e.,
$X^0=\mathbf{0}$. Unless otherwise stated, the Jacobian systems are solved by a two-level restricted additive Schwarz (RAS) preconditioned flexible GMRES (fGMRES) method.
The default restart value of fGMRES is fixed at 50. The ratio of the mesh size between two adjacent levels is 2, two sweeps of RAS preconditioning are applied as the fine-level smoother, and the coarse-level problem is solved approximately by 10 GMRES iterations. A point-block ILU factorization with 3 fill-in levels is used for the subdomain solver. The relative and absolute tolerances of the global nonlinear solver are $\gamma_r=10^{-12}$ and $\gamma_a=10^{-8}$, respectively. In all experiments, the number of MPI processes equals the number of
subdomains $n_p$, with one subdomain assigned to each process. The local
extraction models are trained in parallel on GPUs using the Adam optimizer.

We first introduce the notation used in this section. For the proposed algorithm (Algorithm~\ref{alg:nepinl}), $NI_{col}$ denotes the number of global Newton iterations performed in Phase 1 for collecting residual snapshots used in training, and $NI_{res}$ denotes the number of Newton iterations required in the restart phase. The global Newton iteration count is reported in the form
$NI_{res}(NI_{tot})$, with 
$NI_{tot}=NI_{col}+NI_{res}$,
where $NI_{tot}$ is the total number of global Newton iterations. For the original NE preconditioner (Algorithm~\ref{alg:NE}), only a single number is reported for the global Newton iteration count, since no separate snapshot-collection and restart phases are involved.
$NI_{sub}$ denotes the number of inner subspace Newton iterations. $LI_{sub}$ denotes the average number of fGMRES iterations per Newton step in the subspace nonlinear solve; the same reporting convention is used for the linear iterations in the restart phase. The timings $T_{col}$, $T_{tra}$, $T_{pre}$, and $T_{res}$ represent the compute time in seconds of snapshot collection, model training, nonlinear preconditioning, and the restarted Newton solve, respectively. The total computational time is denoted by $T_{tot}$, with
$$
T_{tot}=T_{col}+T_{tra}+T_{pre}+T_{res},
$$
where the small discrepancies are due to rounding.

\begin{figure*}[!htbp]
  \centering
  \begin{subfigure}{0.28\textwidth}
    \includegraphics[width=\linewidth]{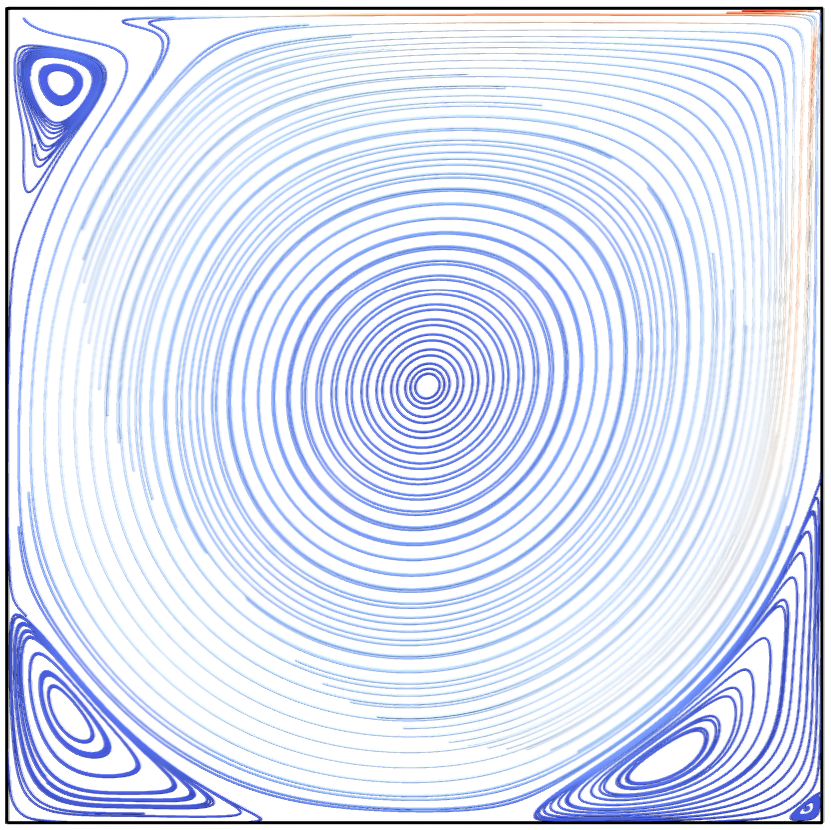}
    \caption{$Re=5{,}000$}
  \end{subfigure}
  \hfill
  \begin{subfigure}{0.28\textwidth}
    \includegraphics[width=\linewidth]{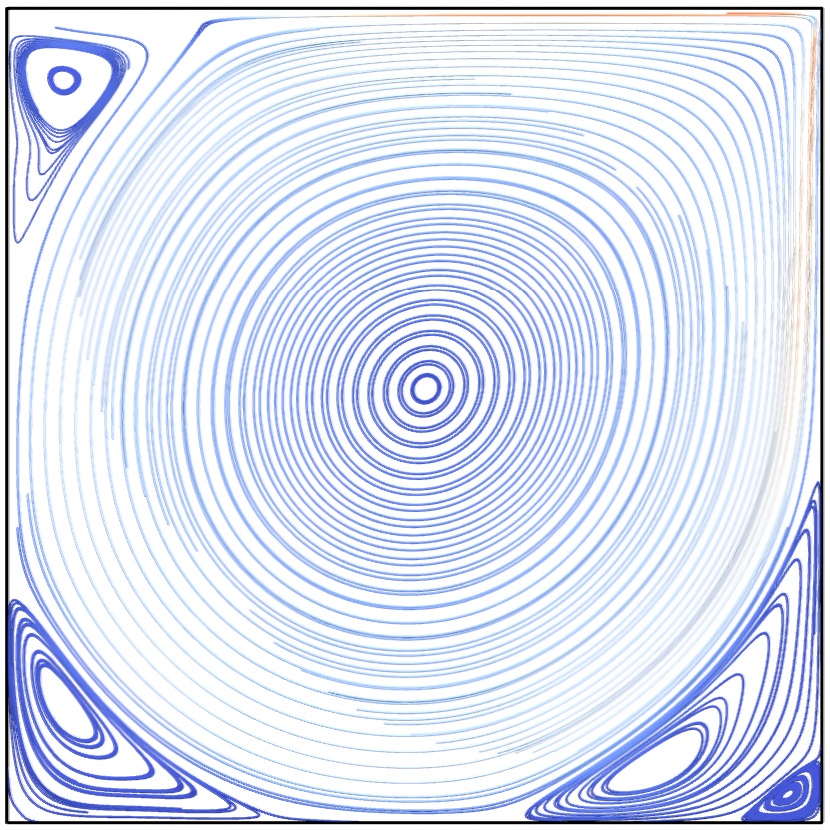}
    \caption{$Re=7{,}500$}
  \end{subfigure}
  \hfill
  \begin{subfigure}{0.28\textwidth}
    \includegraphics[width=\linewidth]{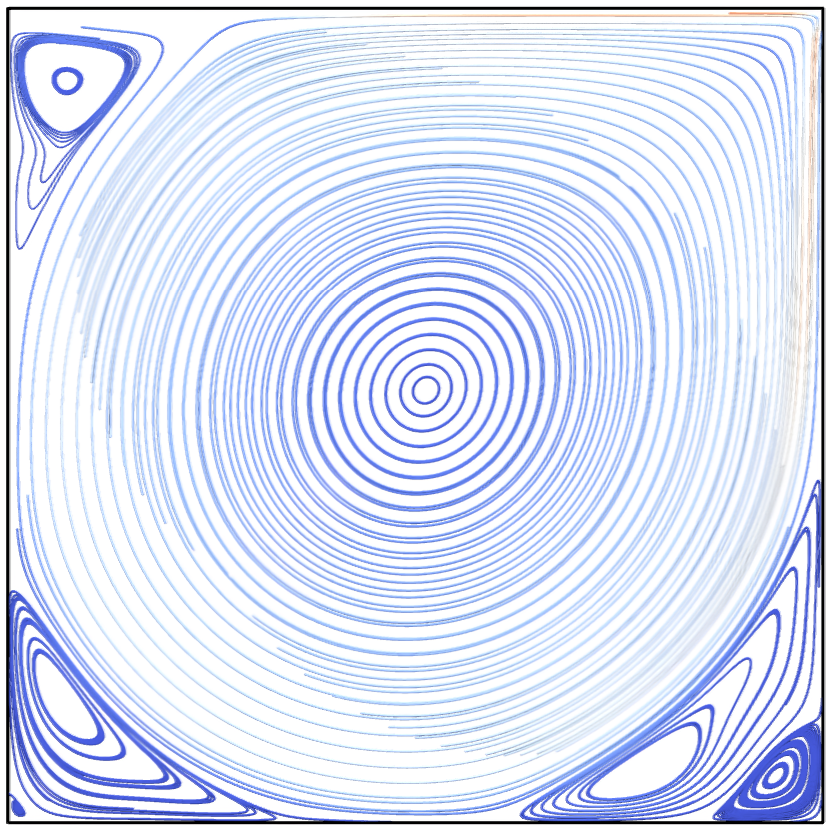}
    \caption{$Re=10{,}000$}
  \end{subfigure}
  \caption{Streamlines of the lid-driven cavity flow obtained by using the  DD-OLE-NE preconditioned inexact Newton method on the $257\times 257$ mesh.}
  \label{fig:streamlines}
\end{figure*}

\begin{figure*}[!htbp]
  \centering
  \begin{subfigure}{0.48\textwidth}
    \includegraphics[width=\linewidth]{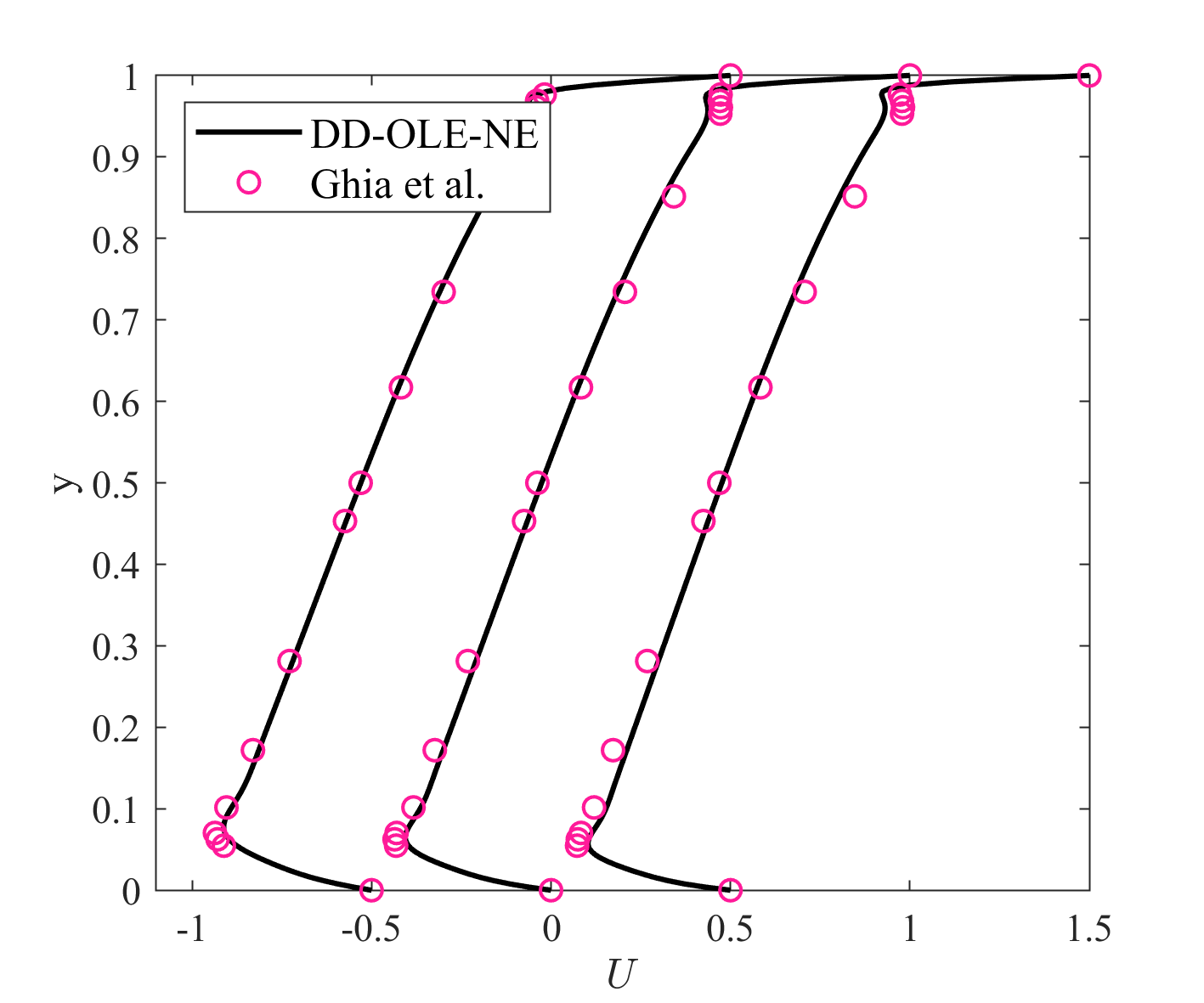}
    \caption{$u$ along the vertical centerline}
  \end{subfigure}
  \hfill
  \begin{subfigure}{0.48\textwidth}
    \includegraphics[width=\linewidth]{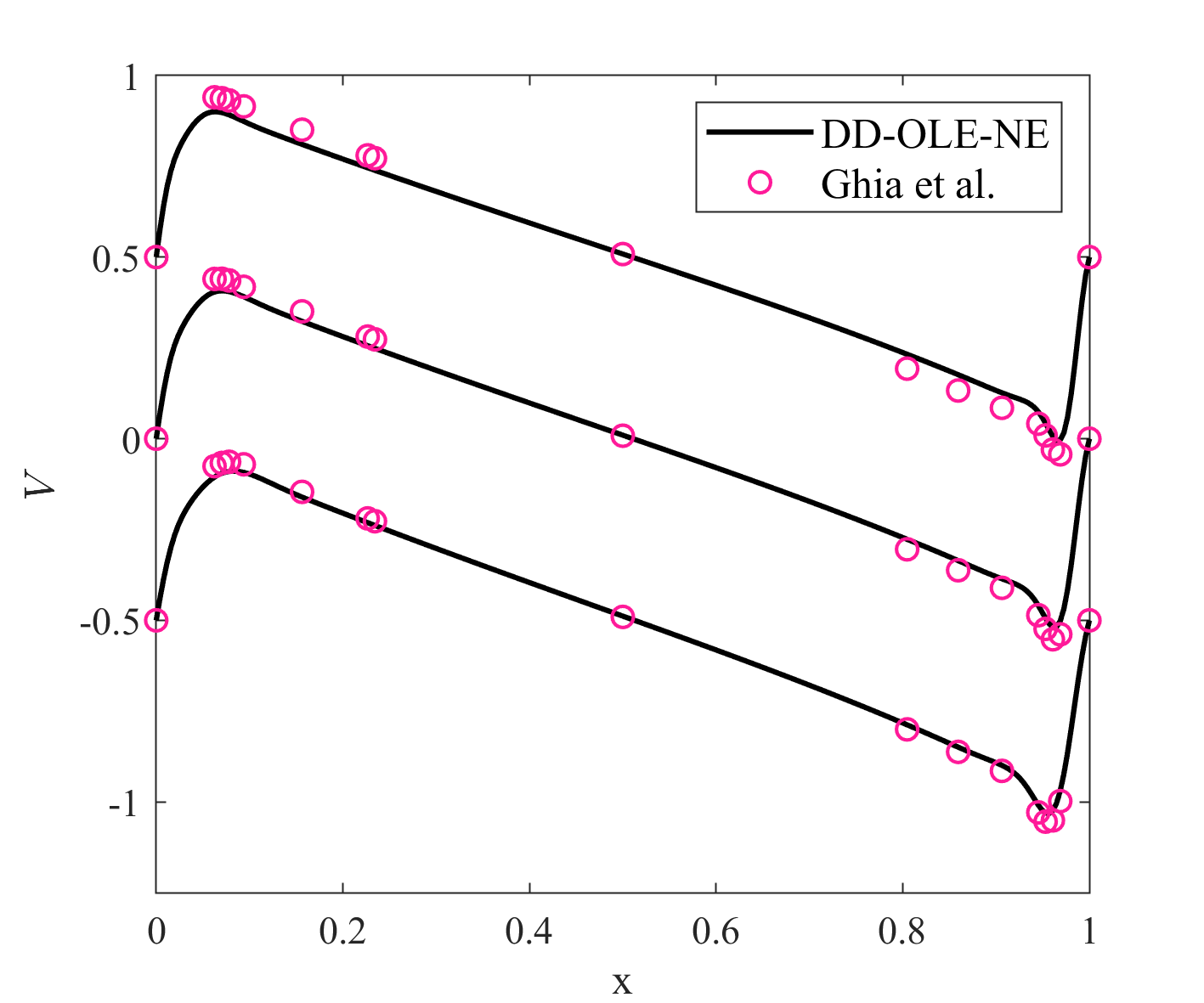}
    \caption{$v$ along the horizontal centerline}
  \end{subfigure}
  \caption{Centerline velocity profiles compared with the results of Ghia et~al.~\cite{Ghia1982} at $Re=5{,}000$ (left/bottom), $7{,}500$ (middle) and $10{,}000$ (right/top).}
  \label{fig:centerline}
\end{figure*}

Figure~\ref{fig:streamlines} presents the streamlines of the lid-driven cavity flows at the three Reynolds numbers considered. Figure~\ref{fig:centerline}  compares the velocity profiles along the centerlines, including the horizontal velocity $u$ along the vertical centerline and the vertical velocity $v$ along the horizontal centerline, against the benchmark data of Ghia et~al.~\cite{Ghia1982}. The computed results agree well with the benchmark profiles for all three Reynolds numbers.

\subsection{Comparison between DD-OLE-NE and baseline NE}\label{subsec:main_comparison}

In this subsection, we compare three variants of the proposed method, namely DD-OLE-NE(PCA), DD-OLE-NE(MLP), and DD-OLE-NE(CNN), with the original NE method in Algorithm~\ref{alg:NE}, which is used as the baseline. To ensure a fair comparison, the baseline NE preconditioner is activated at the same intervention step $k_{re}$ as the proposed DD-OLE-NE methods. In addition, the threshold parameter $\rho$ and the inner nonlinear relative tolerance $\gamma_r^s$ are kept identical for all methods. Unlike the proposed methods, however, the baseline NE method does not require a stagnation period for residual-snapshot collection or online training.

\begin{figure*}[!hbtp]
  \centering
  \includegraphics[width=0.9\linewidth]{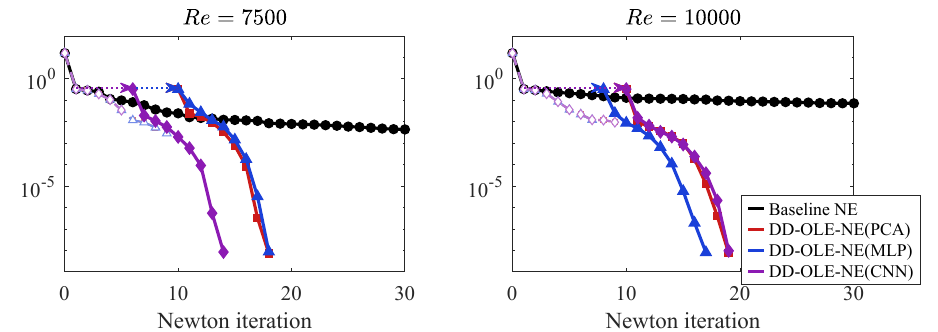}
  \caption{Nonlinear residual histories for baseline NE, DD-OLE-NE(PCA), DD-OLE-NE(MLP), and DD-OLE-NE(CNN) at $Re=7{,}500$ and $10{,}000$. The open markers denote the initial IN phase before nonlinear preconditioning is applied, while the filled markers denote the iterations after applying the preconditioner from the restart iterate $X^{k_{re}}$. The dotted arrows indicate the restart from $X^{k_{re}}$. The corresponding parameters are provided in Table~\ref{tab:main_comparison}.}
  \label{fig:residual_history}
\end{figure*}

Figure~\ref{fig:residual_history} presents the nonlinear residual histories for the lid-driven cavity flow problem at $Re=7{,}500$ and $Re=10{,}000$ on a $257\times257$ mesh. The corresponding choices of $\rho$ and $\gamma_r^s$ are reported in Table~\ref{tab:main_comparison}. For the baseline NE method, the nonlinear residual stagnates around $10^{-2}$, and the iteration fails to reach convergence for both Reynolds numbers. By contrast, the proposed DD-OLE-NE variants learn the dominant residual structures online during the stagnation phase and apply nonlinear elimination at the restart iterate $X^{k_{re}}$. Once the elimination step is activated, the residual decreases rapidly, and the global Newton iteration converges successfully without further stagnation.

Table~\ref{tab:main_comparison} reports the iteration counts and computational performance of the three DD-OLE-NE variants for the two Reynolds numbers. The baseline NE method stagnates and fails to converge under the prescribed values of $\rho$ and $\gamma_r^s$, and is therefore omitted from the table. In contrast, all three proposed variants converge successfully in both cases, demonstrating that extracting the dominant residual structure for the definition of the bad subset can significantly improve the robustness of NE preconditioning.
After the nonlinear elimination restart, the restarted Newton iteration converges rapidly for all tested cases, requiring only $9$--$10$ additional global Newton steps. The total number of global Newton iterations remains moderate, ranging from $14$ to $19$, which indicates relatively stable nonlinear convergence. The three extraction strategies exhibit different advantages. For $Re=7{,}500$, the CNN-based method gives the shortest total computational time and the smallest total number of Newton iterations. For $Re=10{,}000$, the MLP-based method achieves the lowest total time, whereas the PCA-based method introduces the smallest learning overhead. 

The table also shows that the online extraction and nonlinear elimination stages introduce only a modest additional cost compared with the full nonlinear solve. Most of the computational effort is still spent on residual collection and on the restarted global Newton iteration. In particular, for $Re=10{,}000$, the increase in total cost is mainly reflected in the restarted linear solves rather than in the number of nonlinear iterations. This observation suggests that the proposed DD-OLE-NE methods continue to provide a robust nonlinear convergence mechanism when the nonlinear problem becomes more difficult and the associated linearized systems become more expensive to solve. We note that the network architectures used in this and the subsequent subsections are not tuned for each individual test case; their influence on the solver performance is examined in Section~\ref{subsec:nn_params}.

\begin{table*}[!hbtp]\small
   \caption{\textit{The number of iterations and computational times obtained by DD-OLE-NE(PCA), DD-OLE-NE(MLP), and DD-OLE-NE(CNN) for the cavity flow problem with $n_p=4$. The mesh size is $257 \times 257$, $k_s=1$ and $k_{re}=1$. The MLP architecture is $n_i$-$16$-$8$-$16$-$n_i$ where $n_i$ denotes the local residual dimension, trained for $500$ epochs, and the CNN architecture is $3$-$8$-$32$-$8$-$3$ with kernel sizes $7$-$4$-$4$-$7$, trained for $150$ epochs.}}
  \label{tab:main_comparison}
  \centering
  \begin{tabular}{rcccclllllllll}
  \hline
  $Re$ & $d$ & $k_t$ & $\rho$ & $\gamma_r^s$
  & $T_{col}$ & $T_{tra}$ & $NI_{sub}$ & $LI_{sub}$ & $T_{pre}$
  & $NI_{res}(NI_{tot})$ & $LI_{res}$ & $T_{res}$ & $T_{tot}$ \\
  \hline
  \multicolumn{14}{c}{DD-OLE-NE(PCA)}\\
  \hline
  $7{,}500$ & 5 & 9 & 0.35 & 0.10 & 9.7  & 0.04 & 5  & 2.4 & 1.6 & 9(18)  & 5.6  & 4.0  & 15.3 \\
  $10{,}000$            & 8 & 9 & 0.20 & 0.01 & 16.7 & 0.07 & 19 & 2.9 & 6.7 & 10(19) & 22.9 & 15.0 & 38.5 \\
  \hline
  \multicolumn{14}{c}{DD-OLE-NE(MLP)}\\
  \hline
  $7{,}500$ & -- & 9 & 0.35 & 0.10 & 9.7  & 1.7 & 5  & 2.4 & 1.7 & 9(18)  & 6.6  & 5.7  & 18.7 \\
  $10{,}000$            & -- & 7 & 0.20 & 0.01 & 8.5  & 1.8 & 13 & 3.2 & 5.0 & 10(17) & 24.1 & 17.3 & 32.6 \\
  \hline
  \multicolumn{14}{c}{DD-OLE-NE(CNN)}\\
  \hline
  $7{,}500$ & -- & 5 & 0.35 & 0.10 & 4.9  & 1.8 & 4  & 3.0 & 1.5 & 9(14)  & 6.2  & 4.9  & 13.0 \\
  $10{,}000$            & -- & 9 & 0.20 & 0.01 & 16.7 & 2.6 & 19 & 3.1 & 7.2 & 10(19) & 23.3 & 16.6 & 43.0 \\
  \hline
  \end{tabular}
  \end{table*}

\begin{figure*}[htbp]
    \centering
    \begin{subfigure}[t]{0.20\linewidth}
        \centering
        \includegraphics[width=\linewidth]{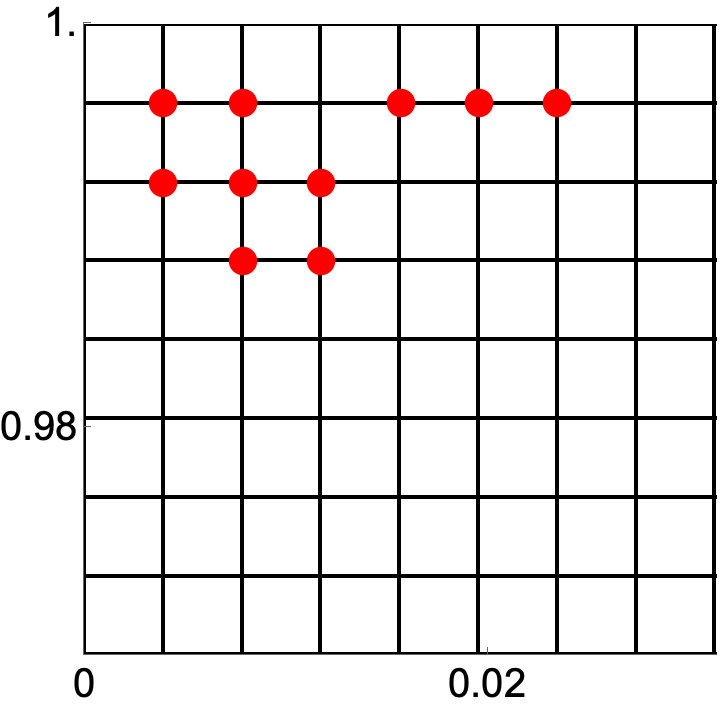}
        \caption{$N_b$}
        \label{fig:rb-nonn}
    \end{subfigure}
    \hfill
    \begin{subfigure}[t]{0.35\linewidth}
        \centering
        \includegraphics[width=\linewidth]{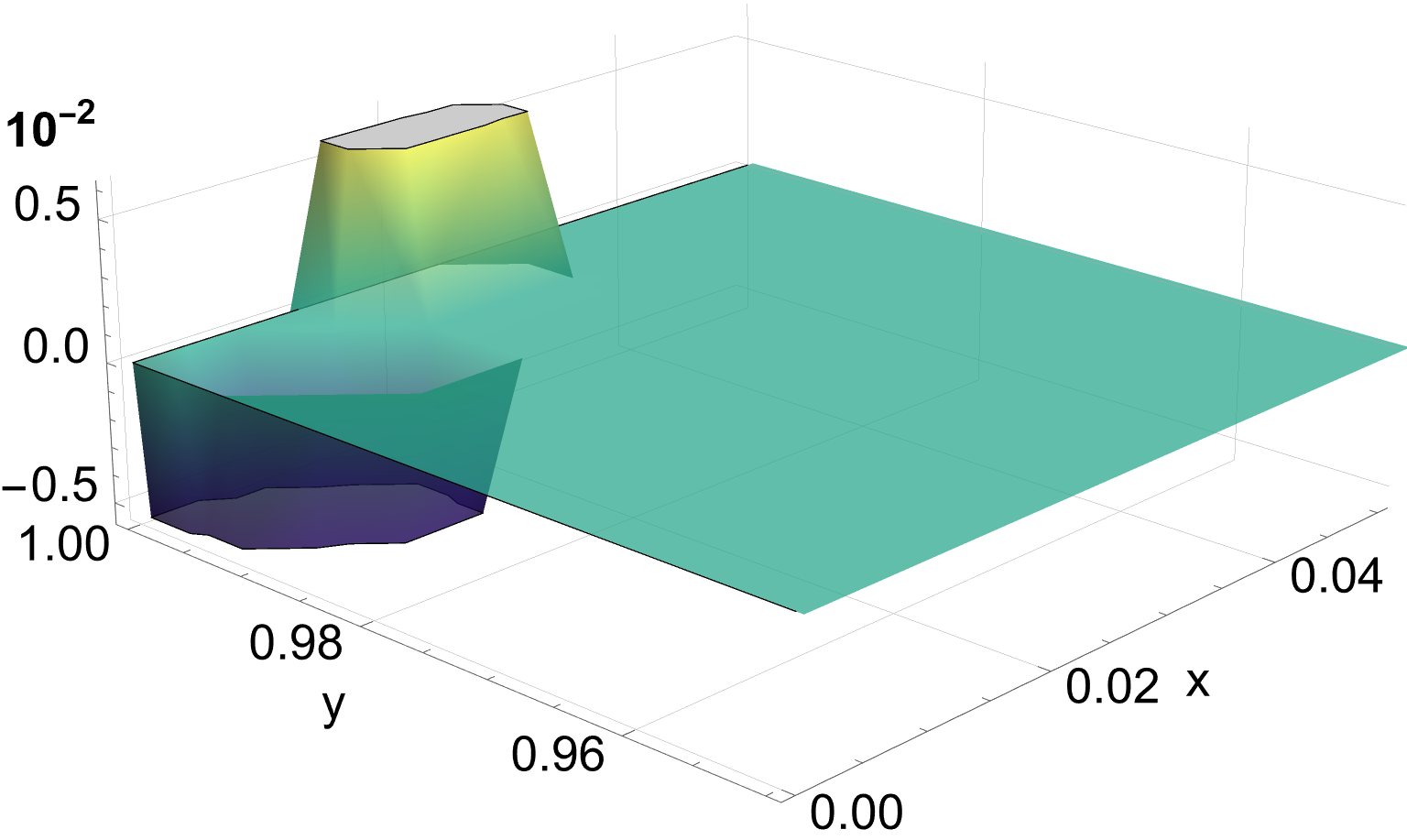}
        \caption{$\bigl(R_bF(Y^0)\bigr)^\omega$}
        \label{fig:nonn-y0}
    \end{subfigure}
    \hfill
    \begin{subfigure}[t]{0.35\linewidth}
        \centering
        \includegraphics[width=\linewidth]{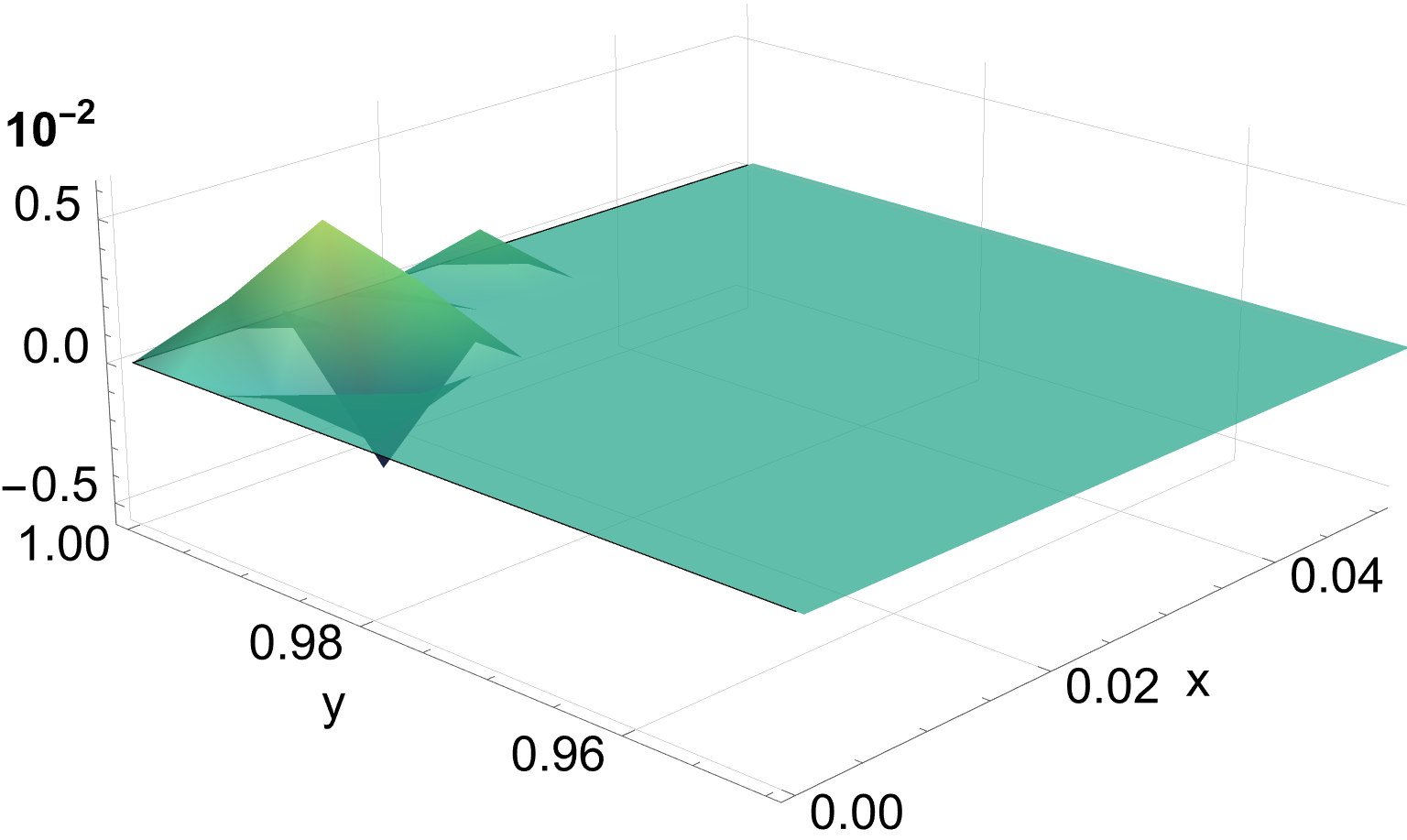}
        \caption{$\bigl(R_bF(Y^\ast)\bigr)^\omega$}
        \label{fig:nonn-ystar}
    \end{subfigure}
    \vspace{0.8cm}
    \begin{subfigure}[t]{0.20\linewidth}
        \centering
        \includegraphics[width=\linewidth]{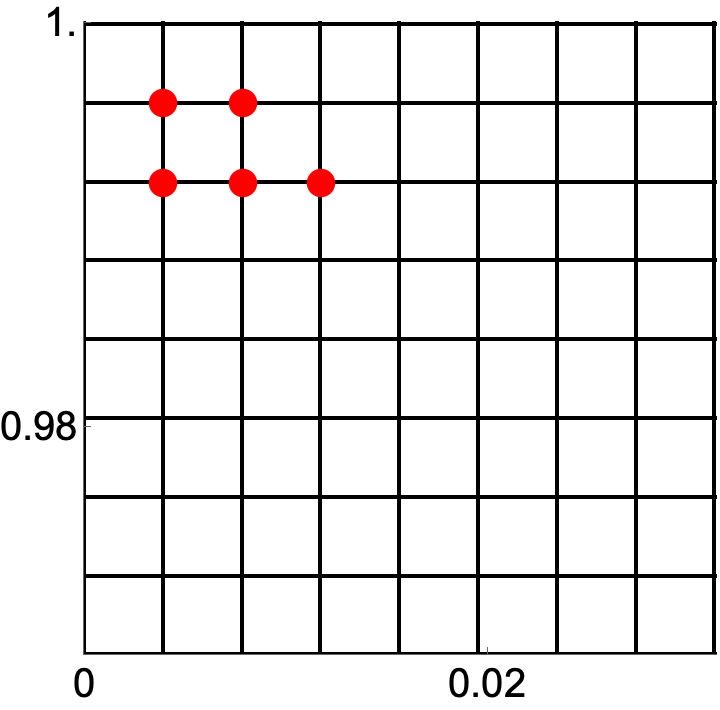}
        \caption{$\mathcal{N}_b$}
        \label{fig:rb-nn}
    \end{subfigure}
    \hfill
    \begin{subfigure}[t]{0.35\linewidth}
        \centering
        \includegraphics[width=\linewidth]{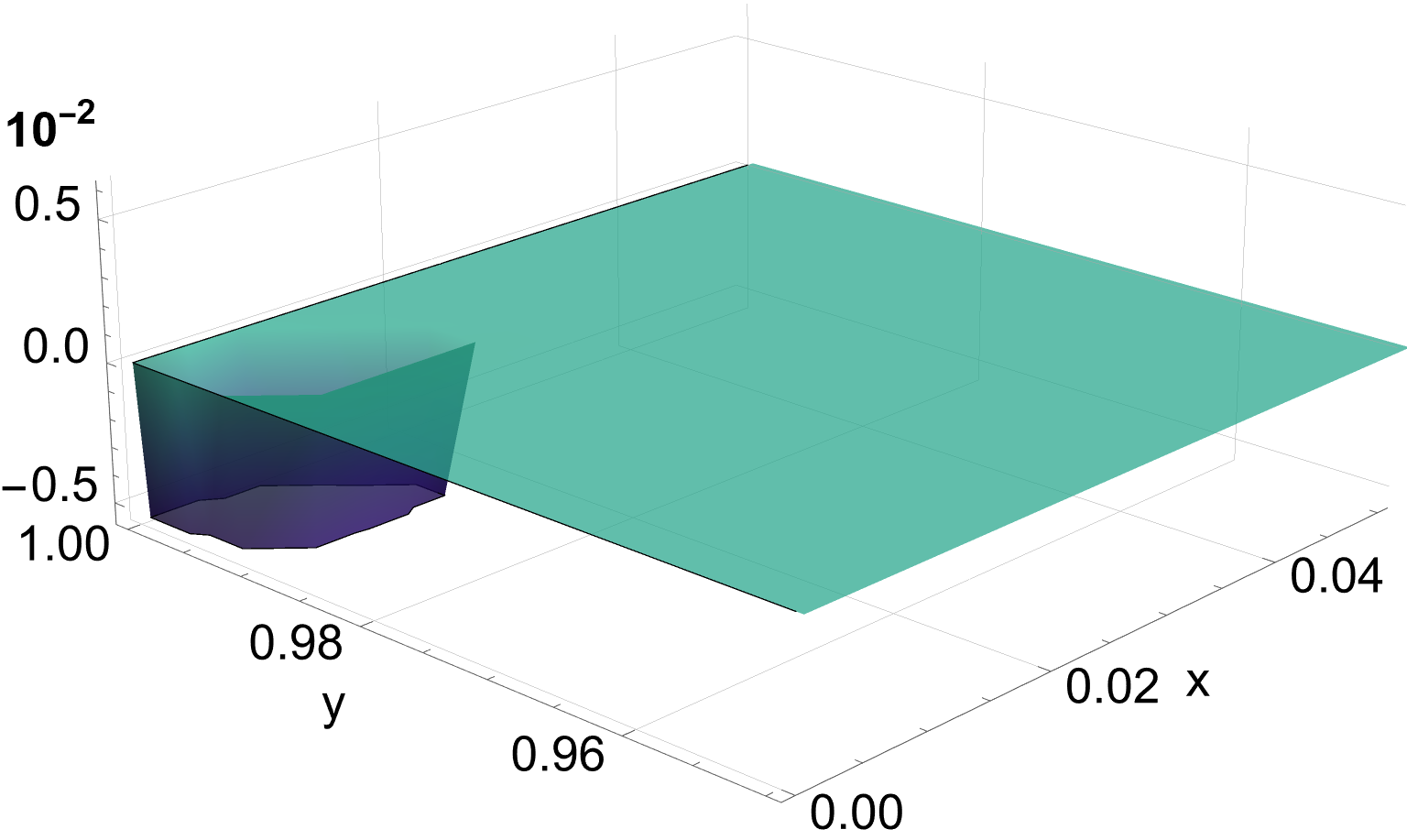}
        \caption{$\bigl(\mathcal{R}_b\mathcal{F}(Y^0)\bigr)^\omega$}
        \label{fig:nn-y0}
    \end{subfigure}
    \hfill
    \begin{subfigure}[t]{0.35\linewidth}
        \centering
        \includegraphics[width=\linewidth]{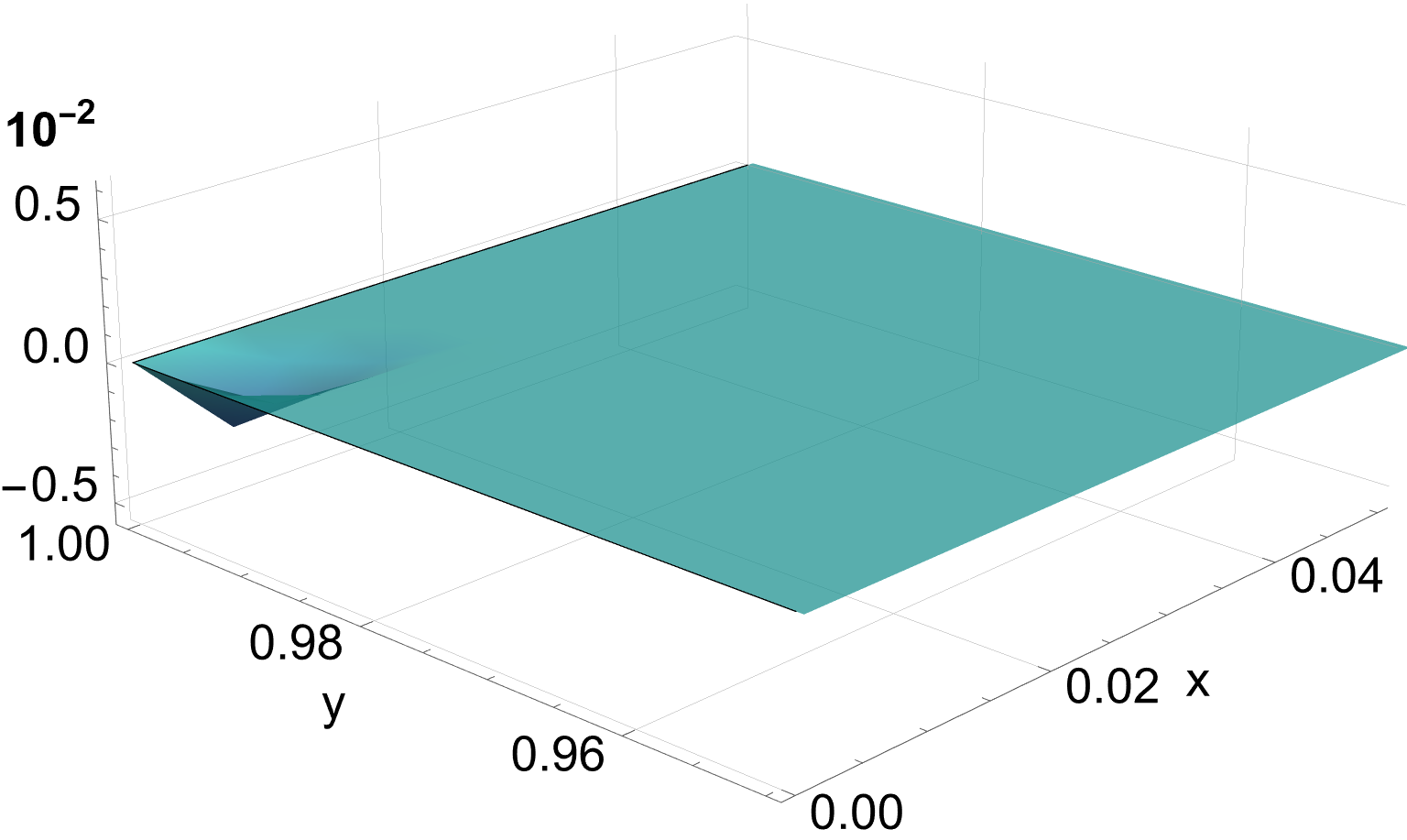}
        \caption{$\bigl(\mathcal{R}_b\mathcal{F}(Y^\ast)\bigr)^\omega$}
        \label{fig:nn-ystar}
    \end{subfigure}
    \caption{Comparison of the selected bad subset and the corresponding restricted vorticity residual surfaces at $Re=10{,}000$. Panels (a)--(c) correspond to the baseline NE method, where the bad subset is selected from the raw residual. Panels (d)--(f) correspond to DD-OLE-NE(MLP), where the bad subset is selected from the learned residual representation. Panels (a) and (d) show the selected bad subset points. Panels (b) and (e) show the restricted residuals at the initial subspace iterate $Y^0$, while panels (c) and (f) show the restricted residuals after the subspace solve, denoted by $Y^\ast $.}

    \label{fig:six-plots}
\end{figure*}

  To further illustrate the effect of the learned residual representation, Figure~\ref{fig:six-plots} compares the bad subset selection and the corresponding restricted vorticity residual surfaces for the baseline NE method and DD-OLE-NE(MLP) at $Re=10{,}000$. When the bad subset is selected from the raw residual, the selected points are more spread out and the restricted residual remains relatively large after the subspace solve. In contrast, the learned residual representation identifies a more compact bad subset associated with the dominant residual structure. After the corresponding subspace solve, the restricted residual is substantially reduced. This comparison provides a local explanation of why applying nonlinear elimination to the learned residual representation leads to more robust convergence than applying NE directly to the raw residual. 
  
\subsection{Effect of neural network extraction}\label{subsec:ablation}

It is well known that neural networks are capable of capturing dominant low-frequency structures from high-dimensional data \cite{rahaman2019spectral}. To determine whether the improvement obtained with neural-network-based extraction arises from learning the dominant residual structures, rather than merely from applying a fixed smoothing operation, we vary only the extraction operator used during the nonlinear elimination phase. All other components of the solver and all algorithmic parameters are kept unchanged, with $\rho=0.35$ and $\gamma_r^s=0.05$. As reported in Table~\ref{tab:ablation_recon}, we compare fixed, data-independent Gaussian smoothing operators with data-adaptive nonlinear autoencoders based on MLP and CNN architectures. The baseline NE preconditioner, which performs nonlinear elimination directly from the raw residual without an extraction step, is also included for reference.

For the Gaussian control, the smoothed residual for each residual component $F^m$ ($m=u,v,\omega$) at an interior grid point $(i,j)$ is defined by
\begin{equation}
  \widetilde{F}^m_{ij}
  =
  \frac{
  \displaystyle
  \sum_{p=-r}^{r}
  \sum_{q=-r}^{r}
  \exp\!\left(-\frac{p^2+q^2}{2\sigma_g^2}\right)
  F^m_{i+p,j+q}
  }{
  \displaystyle
  \sum_{p=-r}^{r}
  \sum_{q=-r}^{r}
  \exp\!\left(-\frac{p^2+q^2}{2\sigma_g^2}\right)
  }.
\end{equation}
Near the boundary, the same formula is applied using only the available neighboring grid points and the weights are renormalized accordingly. Here $\sigma_g$ controls the width of the Gaussian kernel and $r$ is the stencil radius, corresponding to a local $(2r+1)\times(2r+1)$ stencil in the interior. 

In this Gaussian approach, the smoothed residual $\widetilde{F}$ plays the same role as the reconstructed residual $\mathcal{F}$ in the learned DD-OLE-NE variants. It is used to identify the bad subset and is also used in the corresponding restricted nonlinear elimination correction solve. Thus, this ablation tests whether a fixed local smoothing operator can replace the learned extraction.

The results show that only the learned MLP and CNN extractors achieve convergence at both reported Reynolds numbers, whereas the baseline NE preconditioner and both Gaussian smoothing operators fail to converge. This comparison demonstrates that the improvement cannot be attributed solely to replacing the raw residual with a locally smoothed residual. In particular, Gaussian smoothing is a fixed linear extraction procedure whose weights are prescribed independently of the nonlinear iteration, while the MLP and CNN extractors are trained on residual snapshots generated during the stagnation phase and therefore adapt to the residual structures encountered in the current solve. Since neither Gaussian configuration restores convergence, the enhanced robustness is attributed to the ability of the learned extractors to identify problem-dependent dominant residual structures, rather than to fixed residual smoothing alone.

\begin{table*}[!hbtp]\small
\caption{\textit{The number of iterations and computational times obtained using different extraction operators for the cavity flow problem with $n_p=4$. The mesh size is $257\times 257$, $\rho=0.35$, $\gamma_r^s=0.05$, $k_s=1$, $k_t=9$, and $k_{re}=1$. The MLP architecture is $n_i$-$16$-$8$-$16$-$n_i$, trained for $500$ epochs, and the CNN architecture is $3$-$25$-$75$-$25$-$3$ with kernel sizes $5$-$5$-$5$-$5$, trained for $150$ epochs. ``$\times$'' means that the solver does not converge.}}
\label{tab:ablation_recon}
\centering
\setlength{\tabcolsep}{3.5pt}
\begin{tabular}{l ccc ccc}
\hline
 & \multicolumn{3}{c}{$Re=7{,}500$} & \multicolumn{3}{c}{$Re=10{,}000$} \\
\cline{2-4}\cline{5-7}
Methods 
& $NI_{res}(NI_{tot})$ & $LI_{res}$ & $T_{tot}$ 
& $NI_{res}(NI_{tot})$ & $LI_{res}$ & $T_{tot}$ \\
\hline
Baseline NE                 
& \multicolumn{3}{c}{$\times$} 
& \multicolumn{3}{c}{$\times$} \\
\hline
Gaussian smoothing ($\sigma_g=2$, $r=2$)       
& \multicolumn{3}{c}{$\times$} 
& \multicolumn{3}{c}{$\times$} \\
Gaussian  smoothing ($\sigma_g=3$, $r=5$)       
& \multicolumn{3}{c}{$\times$} 
& \multicolumn{3}{c}{$\times$} \\
\hline
DD-OLE-NE(MLP)                       
& 9(18) & 6.3 & 10.1  
& 11(20) & 15.6 & 26.8 \\ 
\hline
DD-OLE-NE(CNN)                       
& 9(18) & 6.8 & 12.1 
& 12(21) & 16.9 & 31.5 \\
\hline
\end{tabular}
\end{table*}

\subsection{Effect of domain  decomposition}\label{subsec:scalability}

We further examine the effect of domain decomposition on DD-OLE-NE(CNN) for the cavity flow problem at $Re=5{,}000$ on a finer mesh, $513\times513$. For this test, we use a three-level hierarchy in the RAS-preconditioned fGMRES solver and set the fGMRES restart value to 100. All other solver parameters are kept unchanged as in the default setting. As shown in Table~\ref{tab:cnn_scaling_re5000_513}, the performance of DD-OLE-NE(CNN) depends strongly on the subdomain partition. When only two subdomains are used, the method fails to converge, indicating that an overly coarse decomposition does not sufficiently localize the dominant nonlinear structures. In this case, each local CNN must represent a relatively large and heterogeneous residual field, which makes the extraction of the relevant nonlinear features more difficult.
Once the domain is divided into four or more subdomains, robust convergence is recovered. As the number of subdomains increases from $4$ to $16$, the restarted Newton iteration count and the average GMRES iteration count remain relatively stable. The total computational time is consequently reduced from $45.7\,s$ to $12.7\,s$.

These results demonstrate two complementary advantages of a finer domain decomposition. First, the local residual snapshots become lower-dimensional and more spatially localized, enabling each subdomain extraction model to identify the dominant nonlinear structure more accurately and making the associated subspace nonlinear problem easier to solve. Second, the local model construction and nonlinear elimination operations can be carried out independently and in parallel across the subdomains. Therefore, within the range tested, using more subdomains improves not only the available parallelism but also the effectiveness of the local extraction models, whereas using too few subdomains may lead to a loss of convergence.

\begin{table}[!hbtp]\small
  \caption{\textit{Effect of subdomain partition on DD-OLE-NE(CNN) for the cavity flow problem at $Re=5{,}000$. The mesh size is $513\times 513$, $\rho=0.10$, $\gamma_r^s=0.01$, $k_{re}=1$, $k_s=1$, and $k_t=9$. A three-level RAS-preconditioned fGMRES solver with restart value 100 is used as the linear solver. The CNN architecture is $3$-$8$-$32$-$8$-$3$ with kernel sizes $7$-$4$-$4$-$7$, trained for 150 epochs. ``$\times$'' means that the case does not converge.}}
  \label{tab:cnn_scaling_re5000_513}
  \centering
  \begin{tabular}{ccccc c}
  \hline
  $n_p$ & $NI_{sub}$ & $LI_{sub}$ & $NI_{res}$ & $LI_{res}$ & $T_{tot}$ \\
  \hline
  2 $(2\times 1)$ & \multicolumn{5}{c}{$\times$} \\
  4 $(2\times 2)$ & 13 & 3.8 & 11 & 5.3 & 45.7 \\
  8 $(4\times 2)$ & 8  & 1.6 & 9  & 6.0 & 17.1 \\
  16 $(4\times 4)$ & 4  & 1.0 & 10 & 5.5 & 12.7 \\
  \hline
  \end{tabular}
\end{table}

\subsection{Effect of the restart point and snapshot collection window}\label{subsec:kre_kt}

We study the influence of the restart point and the snapshot collection window on the performance of DD-OLE-NE(CNN). In particular, $k_t$ determines the end of the tolerable stagnation period, which is used for collecting residual snapshots, while $k_{re}$ specifies the Newton iterate from which the nonlinear elimination restart is performed. The results for the more difficult problem $Re=10{,}000$ are reported in Table~\ref{tab:kre_kts}.

The results show that an earlier restart is generally more robust. When $k_{re}=1$, the solver converges for all tested values of $k_t$. In contrast, larger restart indices lead to less stable behavior, and no convergence is obtained for $k_{re}=7$. This indicates that restarting too late from a strongly stagnated Newton trajectory may reduce the effectiveness of the learned nonlinear correction.

The choice of $k_t$ also affects the overall efficiency. Increasing $k_t$ provides more residual snapshots for training, but it does not necessarily improve the restarted solve. Instead, a longer snapshot-collection period may increase the total computational cost. These results suggest that a relatively short stagnation window is sufficient to extract useful residual structures, and that early restart is preferable in the present test.

\begin{table*}[!hbtp]
\caption{\textit{The impact of $k_{re}$ and $k_{t}$ on the performance of DD-OLE-NE(CNN) for the cavity flow problem with $Re=10{,}000$ and $n_p=4$. The mesh size is $257\times 257$, $k_{s}=1$, $\rho=0.35$ and $\gamma_r^s=0.05$. The CNN architecture is $3$-$8$-$16$-$8$-$3$ with kernel sizes $9$-$7$-$7$-$9$, trained for $300$ epochs. ``$\times$" means the solver does not converge.}} 
\label{tab:kre_kts}
\centering
\scriptsize
\setlength{\tabcolsep}{3pt}
\begin{tabular}{c|ccccccccc}
\hline
$k_{re}$ & \multicolumn{3}{c}{1}      & \multicolumn{3}{c}{4}   & \multicolumn{3}{c}{7}    \\ \hline
$k_{t}$ & $NI_{sub}$ & $NI_{res}(NI_{tot})$ & $T_{tot}$ & $NI_{sub}$ & $NI_{res}(NI_{tot})$ & $T_{tot}$ & $NI_{sub}$ & $NI_{res}(NI_{tot})$  & $T_{tot}$ \\ \hline
3        & 6      & 11(14) & 18.9  & \multicolumn{3}{c}{--} & \multicolumn{3}{c}{--}  \\
5        & 7      & 10(15) & 19.0  & 6      & 12(17) & 15.9  & \multicolumn{3}{c}{--}  \\
7        & 9      & 11(18) & 23.6  & 7      & 14(21) & 28.4  & \multicolumn{3}{c}{$\times$} \\
9        & 12     & 11(20) & 31.1  & 7      & 17(26) & 52.5  & \multicolumn{3}{c}{$\times$} \\
\hline
\end{tabular}
\end{table*}

\subsection{Sensitivity to the threshold and relative tolerance}\label{subsec:rho_gammas}

\begin{figure*}[!htbp]
\centering
\begin{subfigure}{0.48\textwidth}
\includegraphics[width=\linewidth]{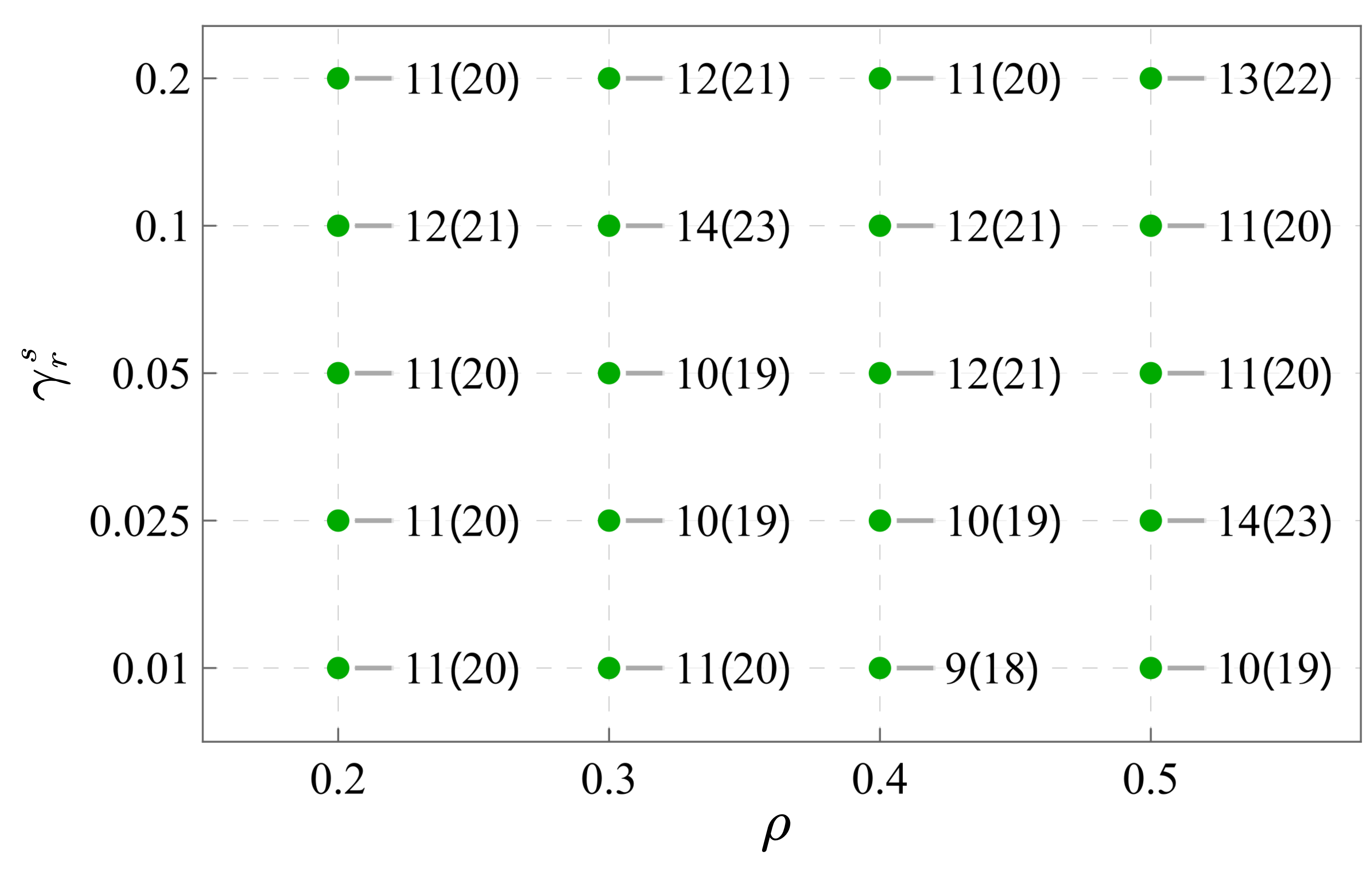}
\caption{DD-OLE-NE(CNN): converges in all $20$ cases}
\end{subfigure}
\hfill
\begin{subfigure}{0.48\textwidth}
\includegraphics[width=\linewidth]{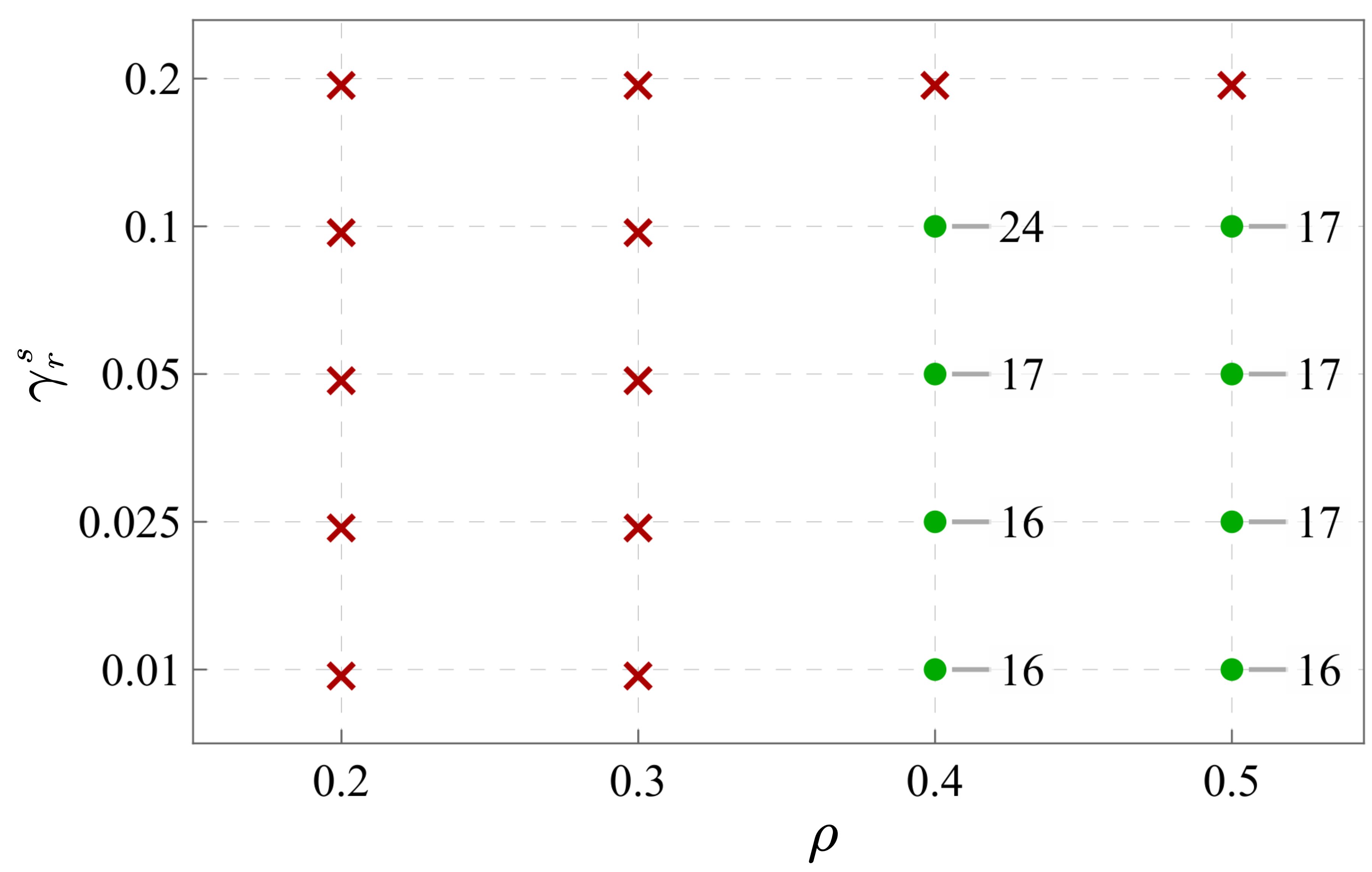}
\caption{NE: converges in $8$ of $20$ cases}
\end{subfigure}
\caption{Parameter sensitivity at $Re=10{,}000$ for DD-OLE-NE(CNN) and baseline NE over the threshold $\rho$ and the inner Newton relative tolerance $\gamma_r^s$. Converged runs are annotated with the Newton iteration count, following the $NI_{res}(NI_{tot})$ convention for DD-OLE-NE(CNN) and using a single number for baseline NE. The CNN architecture is $3$-$25$-$75$-$25$-$3$ with kernel sizes $5$-$5$-$5$-$5$ and is trained for $150$ epochs. The symbol ``$\times$'' means that the case does not converge.}
\label{fig:rho_gammas}
\end{figure*}

We assess the sensitivity of DD-OLE-NE(CNN) to the threshold $\rho$ and the inner Newton relative tolerance $\gamma_r^s$ at $Re=10{,}000$. The parameters are varied over $\rho \in \{0.2,0.3,0.4,0.5\}$ and $\gamma_r^s \in \{0.01,0.025,0.05,0.1,0.2\}$. The results are compared with the baseline NE method in Figure~\ref{fig:rho_gammas}.

For DD-OLE-NE(CNN), residual snapshots are collected during the initial IN stagnation window from $k_s=1$ to $k_t=9$. After the CNN extractor is trained, the nonlinear elimination correction is applied by restarting from $X^{k_{re}}$ with $k_{re}=1$. In contrast, the baseline NE does not use the learned extraction or the early restart, and the NE preconditioner is applied directly to the raw residual starting from the ninth Newton iteration.

DD-OLE-NE(CNN) converges for all tested parameter pairs. Moreover, the restarted Newton iteration count remains relatively stable over the whole parameter range, indicating that the learned extraction reduces the sensitivity of the bad-subset selection to the choices of $\rho$ and $\gamma_r^s$. In contrast, the baseline NE method converges only for part of the parameter space. In particular, it fails for all tested tolerances when $\rho=0.2$ or $\rho=0.3$, and also fails for the loosest tolerance $\gamma_r^s=0.2$ when $\rho=0.4$ or $0.5$.

These results show that the main advantage of DD-OLE-NE(CNN) is not merely a reduction in iteration count at a single parameter setting, but a significantly more robust convergence behavior over the parameter space. By applying the bad-subset selection rule to a learned residual representation rather than to the raw residual, DD-OLE-NE(CNN) provides a more stable nonlinear preconditioning strategy.

\subsection{Impact of neural network parameters}\label{subsec:nn_params}

We further examine how the network architecture and training affect the performance of the proposed method. The study includes the training-epoch budget, the convolutional kernel sequence, and the channel sequence. All tests in this subsection are conducted for the cavity flow problem at $Re=10{,}000$ with $n_p=4$ on a $257\times257$ mesh.

Table~\ref{tab:epoch_impact} compares the effect of the training-epoch budget for the MLP and CNN extractors. The CNN-based method converges for all tested epoch numbers, while the MLP-based method requires a sufficiently large training budget to achieve convergence. This indicates that, in this test, the CNN extractor is less sensitive to the number of training epochs. However, increasing the epoch number does not lead to a monotonic reduction in the total computational time, suggesting that excessive training is not necessarily beneficial for the overall solver performance.

Tables~\ref{tab:kernel_size} and~\ref{tab:channels} report the influence of the CNN kernel sizes and channel sequences. The results show that DD-OLE-NE(CNN) remains effective for a wide range of network architectures. In particular, the total compute time is not solely determined by the number of trainable parameters. Moderate-size networks with a bottleneck structure are sufficient to construct an effective residual extractor, whereas simply increasing the network width or parameter count does not guarantee improved performance. In some cases, larger networks may increase the cost of the restarted Newton solve.

Overall, these results suggest that DD-OLE-NE(CNN) is reasonably robust with respect to the neural network design within the tested range. The network should have enough capacity to capture the dominant residual structures, but overly large architectures are unnecessary and may reduce computational efficiency.

\begin{table}[!hbtp]
\caption{\textit{The impact of training epochs on the performance of DD-OLE-NE(MLP) and DD-OLE-NE(CNN) for the cavity flow problem with $Re=10{,}000$ and $n_p=4$. The mesh size is $257\times 257$, $\rho=0.3$, and $\gamma_r^s=0.05$. The MLP architecture is $n_i$-$16$-$8$-$16$-$n_i$, and the CNN architecture is $3$-$25$-$75$-$25$-$3$ with kernel sizes $5$-$5$-$5$-$5$.}}
\label{tab:epoch_impact}
\centering
\begin{tabular}{c|cc|cc}
\hline
\multicolumn{1}{c|}{} & \multicolumn{2}{c|}{MLP} & \multicolumn{2}{c}{CNN} \\ 
\cline{2-3} \cline{4-5} 
$epoch$ & $NI_{sub}$ & $T_{tot}$ & $NI_{sub}$ & $T_{tot}$ \\ 
\hline
150 & \multicolumn{2}{c|}{$\times$} & 12 & 34.8 \\
300 & 14 & 35.3 & 13 & 33.1 \\
500 & 5  & 27.7 & 9  & 41.6 \\
\hline
\end{tabular}
\end{table}

\begin{table*}[!hbtp]
\caption{\textit{Impact of the CNN kernel size sequence on the performance of DD-OLE-NE(CNN) for the cavity flow problem with $Re=10{,}000$ and $n_p=4$. The mesh size is $257\times257$, with $\rho=0.1$, $\gamma_r^s=0.05$, $k_{re}=1$, $k_s=1$, and $k_t=3$. The channel sequence is 3-8-16-8-3 for the first six architectures and 3-8-16-64-16-8-3 for the remaining architectures and the networks are trained for 150 epochs. \#Params denotes the number of parameters for each local neural network.}}
\label{tab:kernel_size}
\centering
\begin{tabular}{cc|ccccc}
\hline
Kernel sizes & \#Params & $NI_{sub}$ & $NI_{res}(NI_{tot})$ & $T_{tra}$ & $T_{res}$ & $T_{tot}$ \\ \hline
5,3,3,5  & 14156  & 6      & 10(13)     & 1.7         & 13.3          & 17.0        \\
5,11,11,5 & 128844   & 6      & 13(16)     & 1.9         & 16.0          & 19.9        \\
7,4,4,7 & 25932    & 7      & 10(13)     & 1.7         & 15.3          & 19.0        \\
7,12,12,7 & 157004   & 8      & 11(14)     & 2.0         & 16.8          & 20.8        \\
9,5,5,9 & 41292     & 7      & 11(14)     & 1.7         & 15.3          & 19.0        \\
9,13,13,9 & 188748   & 6      & 11(14)     & 2.1         & 14.0          & 18.1        \\
3,4,6,6,4,3 & 313484 &   7   &   16(19)  &   1.8 &   21.4    &   25.2    \\
5,5,5,5,5,5 & 235660 &  9   &   13(16)  &   1.8 &   22.3    &   26.1    \\
7,4,3,3,4,7 & 99980 &   5   &   12(15)  &   1.8 &   11.5    &   15.3    \\
7,8,9,9,8,7 & 738956 &  5   &   12(15)  &   1.9 &   18.7    &   22.6    \\
\hline
\end{tabular}
\end{table*}

\begin{table*}[!hbtp]
\caption{\textit{Impact of the CNN channel sequence on the performance of DD-OLE-NE(CNN) for the cavity flow problem with $Re=10{,}000$ and $n_p=4$. The mesh size is $257\times257$, with $\rho=0.1$, $\gamma_r^s=0.05$, $k_{re}=1$, $k_s=1$, and $k_t=3$. The kernel sizes are fixed to 5 in all layers and the networks are trained for $150$ epochs.}}
\label{tab:channels}
\centering
\begin{tabular}{cc|ccccc}
\hline
\multicolumn{1}{c}{Channels} & \#Params & $NI_{sub}$ & $NI_{res}(NI_{tot})$ & $T_{tra}$ & $T_{res}$ & $T_{tot}$ \\ \hline
3,4,8,4,3           & 8876      & 7      & 12(15)     & 1.6         & 16.6          & 20.2         \\
3,8,16,8,3          & 30540      & 7      & 13(16)     & 1.7         & 16.7          & 20.4         \\
3,16,32,16,3        & 112268     & 8      & 16(19)     & 1.7         & 14.8          & 18.5         \\
3,8,32,8,3          & 56204     & 9      & 10(13)     & 1.6         & 16.3          & 19.9         \\
3,8,64,8,3          & 107532     & 5      & 11(14)     & 1.6         & 11.5          & 15.1         \\
3,8,16,64,16,8,3    & 235660     & 9      & 13(16)     & 1.8         & 22.3          & 26.1         \\ \hline
\end{tabular}
\end{table*}

\section{Concluding remarks}
\label{sec:conclusions}
We proposed and studied DD-OLE-NE, an online-learning-enhanced nonlinear elimination preconditioner for inexact Newton methods, targeting nonlinear systems of algebraic equations with strong local nonlinearity. The central idea is to learn the dominant structure of the nonlinear residual online, directly from the residual snapshots that the stagnating global Newton iteration has already generated, and to use the resulting learned extraction to guide the identification of the bad subset eliminated by the nonlinear preconditioner. Because the reconstructed residual exposes the dominant structure of the local high nonlinearity while filtering out undesired perturbations, a more reliable bad subset can be identified with a small threshold and the resulting subspace nonlinear system is more tractable to solve. We considered both a linear extractor based on PCA and nonlinear extractors based on autoencoder neural networks, which are trained independently, so the learning step fits naturally into the existing domain decomposition parallel framework. We evaluated the method on the two-dimensional lid-driven cavity flow at high Reynolds numbers. Compared with the baseline NE preconditioner, DD-OLE-NE converges across a wide range of the NE-related parameters and learning-related parameters. The focus of this paper has been on incompressible flow, but the method is algebraic and is expected to apply to other highly nonlinear problems.

\bibliographystyle{spmpsci}
\bibliography{references}

\end{document}